\documentstyle[12pt]{article}     
\begin{document} 
\title{Stochastic processes on non-Archimedean spaces. 
III. Stochastic processes on totally disconnected topological 
groups.}  
\author{S.V. Ludkovsky}  
\date{17 May 2001
\thanks{Mathematics subject classification
(1991 Revision) 28C20 and 46S10.} }
\maketitle
\par address: Theoretical Department,
Institute of General Physics,
Str. Vavilov 38, Moscow, 119991 GSP-1, Russia.
\begin{abstract} 
Stochastic processes on totally disconnected topological groups
are investigated. In particular they are considered
for diffeomorphism groups and loop groups
of manifolds on non-Archimedean Banach spaces.  
Theorems about a quasi-invariance and a pseudo-differentiability of 
transition measures are proved.  Transition measures are used
for the construction of strongly continuous representations 
including irreducible of these groups.
In addition stochastic processes on general Banach-Lie groups, 
loop monoids, loop spaces
and path spaces of manifolds on Banach spaces over
non-Archimedean local fields also are investigated.
\end{abstract} 
\section{Introduction.}
This part is the continuation of the previous two 
\cite{lunast1,lunast2}, where stochastic processes on
Banach spaces over local fields and stochastic antiderivational
equations on them were investigated.
This part is devoted to stochastic processes on a
totally disconnected topological group which is 
complete, separable and ultrametrizable. 
In particular stochastic processes on 
diffeomorphism groups and loop groups 
of manifolds on Banach spaces over a local field
are considered. These groups were defined and investigated 
in previous articles of the author 
\cite{lud,luumnls,luseamb,lubp2,lutmf99}.
These groups are non-locally compact and 
for them the Campbell-Hausdorff formula is not valid (in 
an open local subgroup). 
In this article topological groups satisfying locally the 
Campbell-Hausdorff formula also are considered.
\par Finite-dimensional Lie groups 
satisfy locally the Campbell-Hausdorff formula.
This is guarantied, if to impose on a locally compact topological
Hausdorff group $G$ two conditions: it is a $C^{\infty }$-manifold
and the following mapping $(f,g)\mapsto f\circ g^{-1}$
from $G\times G$ into $G$ is of class $C^{\infty }$.
But for infinite-dimensional $G$ the Campbell-Hausdorff 
formula does not follow from these conditions.
Frequently topological Hausdorff groups satisfying these 
two conditions also are called Lie groups, though they 
can not have all properties of finite-dimensional 
Lie groups, so that the Lie algebras for them do not play 
the same role as in the finite-dimensional case and therefore
Lie algebras are not so helpful.
If $G$ is a Lie group and its tangent space $T_eG$ is a Banach space,
then it is called a Banach-Lie group, sometimes it is undermined,
that they satisfy the Campbell-Hausdorff formula locally for a Banach-Lie 
algebra $T_eG$. In some papers the Lie group terminology
undermines, that it is finite-dimensional.
It is worthwhile to call Lie groups satisfying the Campbell-Hausdorff
formula locally (in an open local subgroup) by Lie groups in the 
narrow sense;
in the contrary case to call them by Lie groups in the broad sense.
\par In this article also theorems about a quasi-invariance and 
a pseudo-differentiability of transition measures on the totally 
disconnected topological group $G$ relative to the dense subgroup $G'$
are proved.
For measures on Banach spaces over locally compact 
non-Archimedean fields their quasi-invariance and pseudo-differentiability
were investigated in \cite{lu6} (see also
\cite{lutmf99,lubp2} for diffeomorphism and loop groups, but for 
measures not related with stochastic processes).
In each concrete case of  $G$ it its necessary to construct a stochastic 
process and $G'$. 
Below path spaces, loop spaces, 
loop monoids, loop groups and diffeomorphism groups are considered
not only for finite-dimensional, but also for infinite-dimensional 
manifolds.
\par In particular, loop and diffeomorphism groups are important for the 
development of the representation theory of non-locally compact groups. 
Their representation theory has many 
differences with the traditional representation theory of 
locally compact groups and finite-dimensional Lie groups, 
because non-locally 
compact groups have not $C^*$-algebras associated with the Haar 
measures and they have not underlying Lie algebras and relations between 
representations of groups and underlying algebras (see also \cite{lubp}). 
\par In view of the A. Weil 
theorem if a topological Hausdorff group $G$ has a quasi-invariant measure 
relative to the entire $G$, then $G$ is locally compact.
Since loop groups $(L^MN)_{\xi }$ are not locally compact, they can not have 
quasi-invariant measures relative to the entire group, but only relative to 
proper subgroups $G'$ which can be chosen dense in  $(L^MN)_{\xi }$,
where an index $\xi $ indicates on a class of smoothness.
The same is true for diffeomorphism groups.
\par In this article classes of smoothness of the type 
$C^n$ by Schikhof are used. We recall shortly their definition.
Let $\bf K$ be a local field,
that is, a finite algebraic extension of the $p$-adic 
field $\bf Q_p$ for the corresponding prime number $p$ \cite{wei}.
For $b\in \bf R$, 
$0< b <1$, we consider the following mapping:
$$\mbox{ }j_b(\zeta ):= p^{b\times ord_p(\zeta )}\in \bf \Lambda _p$$ 
for $\zeta \ne 0$,  $j_b(0)
:=0$, such that $j_b(*): {\bf K}\to \bf \Lambda _p$, where 
$\bf K\subset \bf C_p$, 
$\bf C_p$ denotes the field of complex numbers with the 
non-Archimedean valuation extending that of $\bf Q_p$, 
$p^{-ord_p(\zeta )}:=|\zeta |_{\bf K},$  
$\bf \Lambda _p$ is a spherically complete
field with a valuation group $ \{ |x|:$ $0\ne x \in {\bf \Lambda _p} \} =
(0,\infty )\subset \bf R$ such that ${\bf C_p}\subset \bf \Lambda _p$
\cite{diar,roo,sch1,wei}. Then we denote $j_1(x):=x$ for each $x\in \bf K$.
Let us consider Banach spaces
$X$ and $Y$ over $\bf K$. Suppose $F: U\to Y$ is a mapping, 
where $U\subset X$ is an open bounded subset.
The mapping $F$ is called differentiable if for each $\zeta \in \bf K$,
$x \in U$ and $h \in X$ with $x+\zeta h \in U$ there exists
a differential such that
$$(1)\mbox{ }DF(x,h):= 
dF(x+\zeta h)/d\zeta  \mid _{\zeta =0 } :=\lim_{\zeta \to 0} 
\{ F(x+\zeta h)-F(x) \}/\zeta $$  and
$DF(x,h)$ is linear by $h$, that is, $DF(x,h)=:F'(x)h$, where $F'(x)$
is a bounded linear operator (a derivative). Let 
$$(2)\mbox{ }\Phi ^1 F(x;h;\zeta )
:=\{ F(x+\zeta h)-F(x)\} /\zeta $$ 
be a partial difference quotient of order $1$
for each $x+\zeta h\in U$, $\zeta h\ne 0$. If $\Phi ^1 F(x;h;\zeta )$ 
has a bounded continuous extension
${\bar \Phi }^1 F$ onto $U\times V\times S$,
where $U$ and $V$ are open neighbourhoods of $x$ and $0$ in $X$,
$U+V\subset U$, $S=B({\bf K},0,1)$, then
$$(3)\mbox{ }\| {\bar \Phi }^1 F(x;h;\zeta ) \| 
:=\sup_{(x \in U, h  \in V, \zeta 
\in S)}\| {\bar \Phi }^1 F(x;h;\zeta ) \| _Y< \infty $$  
and ${\bar \Phi }^1 F(x;h;0)=F'(x)h$. Such $F$ is called
continuously differentiable on $U$. The space of such $F$ is denoted
$C(1,U\to Y)$. Let 
$$(4)\mbox{ }\Phi ^b F(x;h;\zeta ):=(F(x+\zeta h)-F(x))/j_b(\zeta )\in 
Y_{\bf \Lambda _p}$$ be partial difference quotients of order $b$
for $0<b<1$, $x+\zeta h\in U$, $\zeta h\ne 0$, $\Phi ^0F:=F$, where
$Y_{\bf \Lambda _p}$ is a Banach space obtained from $Y$ by 
extension of a scalar field from $\bf K$ to $\bf \Lambda _p$.
By induction using Formulas $(1-4)$
we define partial difference quotients of orders $n+1$ and $n+b$:
$$(5)\mbox{ }\Phi ^{n+1} F(x;h_1,...,
h_{n+1};\zeta _1,...,\zeta _{n+1} ):= $$
$$\{ \Phi ^{n} F(x+\zeta _{n+1}
h_{n+1};h_1,..,h_n;
\zeta _1,...,\zeta _n)- \Phi ^n F(x;h_1,...,h_n;$$
$$\zeta _1,...,\zeta _n) \}/\zeta _{n+1}\mbox{ and }
(\Phi ^{n+b} F)=\Phi ^b( \Phi ^nF)$$ and derivatives
$F^{(n)}=(F^{(n-1)})'.$
Then $C(t,U\to Y)$ is a space of functions $F: U\to Y$ for which
there exist bounded continuous extensions
${\bar \Phi }^vF$ for each $x$ and $x+\zeta _ih_i \in U$
and each $0\le v \le t$, such that each derivative
$F^{(k)}(x): X^k\to Y$ is a continuous $k$-linear operator for each
$x\in U$ and $0<k\le [t],$ where $0\le t<\infty $,
$h_i\in V$ and $\zeta _i\in S,$ 
$[t]=n\le t$ and $\{ t\} =b$ are the integral 
and the fractional parts of $t=n+b$ respectively.
The norm in the Banach space $C(t,U\to Y)$ is the following:
{\large
$$(6)\mbox{ }\| F\|_{C(t,U\to Y)} 
:=sup_{(x, x+\zeta _ih_i\in U; h_i\in V; \zeta _i\in S;
i=1,...,s=[v]+sign \{ v\} ; 0\le v \le t)}$$ 
$$\| ({\bar \Phi }^v F)(
x;h_1,..,h_s;\zeta _1,...,\zeta _s) \|_{Y_{\bf \Lambda _p}},$$ }
where $0\le t\in \bf R,$ $sign (y)=-1$ for $y<0$, $sign (y)=0$ for $y=0$
and $sign (y)=1$ for $y>0$.
\par It is necessary to note that there are quite another groups
with the same name loop groups, but they are infinite-dimensional
Banach-Lie groups of mappings $f: M\to H$ into a finite-dimensional Lie group
$H$ with the pointwise group multiplication of mappings with values in $H$.
The loop groups considered here are geometric loop groups.
\par On the other hand, representation theory of non-locally compact groups
is little developed apart from the case of locally compact groups.
For locally compact groups theory of induced representations
is well developed due to works of Frobenius, Mackey, etc.
(see \cite{fell} and references therein).
But for non-locally compact groups it is very little known.
In particular geometric loop and diffeomoprphism groups have important 
applications in modern physical theories (see \cite{lubp2,lutmf99}
and references therein).
\par Then measures are used for the study of associated unitary 
representations of dense subgroups $G'$. 
\section{Stochastic antiderivational equations and measures
on totally disconnected topological groups.}
\par To avoid misunderstandings we first remind our definitions
from \cite{luseamb,lubp2,lutmf99}.
\par {\bf 2.1. Definitions and Notes. 1.} Let $X$ be a Banach space
over a local field $\bf K$. 
Suppose $M$ is an analytic manifold modelled on $X$ with
an atlas $At(M)$ consisting of disjoint clopen charts $(U_j,\phi _j)$, 
$j\in \Lambda _M$, $\Lambda _M\subset \bf N$. 
That is, $U_j$ and $\phi _j(U_j)$
are clopen in $M$ and $X$ respectively, $\phi _j:U_j\to \phi _j(U_j)$
are homeomorphisms, $\phi _j(U_j)$ are bounded in $X$.
Let $X=c_0(\alpha ,{\bf K})$, where
$$(1)\mbox{ }c_0(\alpha ,{\bf K}):=
\{ x=(x^i:\mbox{ } i \in \alpha )| x^i\in {\bf K},
\mbox{ and for each } \epsilon >0 \mbox{ the set } $$
$$(i: |x^i|> \epsilon )\mbox{ is finite } \} \mbox{ with}$$  
$$(2)\mbox{ } \| x \| :=\sup_i |x^i |<\infty $$ and 
the standard orthonormal base
$(e_i:$ $i\in \alpha )$ \cite{roo}, $\alpha $ is 
considered as an ordinal due to the Kuratowski-Zorn lemma,
$\alpha \ge 1$. Its cardinality is called a dimension
$card (\alpha )=:dim_{\bf K}c_0(\alpha ,{\bf K})$ over $\bf K$. 
\par Then $C(t,M\to Y)$ for $M$
with a finite atlas $At(M)$, $card(\Lambda _M)<\aleph _0$,
denotes a Banach space of functions $f: M\to Y$ with an ultranorm
$$(3)\mbox{ }\| f\| _t=
\sup_{j \in \Lambda _M}\| f|_{U_j} \| _{C(t,U_j\to Y)}< \infty ,$$ where
$Y:=c_0(\beta ,{\bf K})$ is the Banach space over $\bf K$, 
$0\le t\in \bf R$, their restrictions
$f|_{U_j}$ are in $C(t,U_j\to Y)$ for each $j$, $\beta \ge 1$. 
\par  {\bf 2.1.2.} Let $X,$ $Y$ and $M$ be the same as in \S 2.1.1
for a local field $\bf K$. 
We denote by $C_0(t,M\to Y)$ a completion of a subspace of 
cylindrical functions restrictions of which 
on each chart $f|_{U_l}$ are finite $\bf K$-linear combinations of 
functions $\{ {\bar Q}_{\bar m}(x_{\bar m})q_i|_{U_l}:\mbox{ } i\in \beta , 
m \} $ relative to the following norm:
$$(1)\mbox{ }\| f\| _{C_0(t,M\to Y)}:=\sup_{i,m,l}|a(m,f^i|_{U_l})|J_l(t,m),$$
where multipliers $J_l(t,m)$ are defined as follows:
$$(2)\mbox{ }J_l(t,m):= \| {\bar Q}_{\bar m}|_{U_l} \|_{
C(t,\phi _l(U_l)\cap {\bf K^n}\to {\bf K})},$$
$m\in c_0(\alpha ,{\bf Q_p})$ with components $m_i\in \bf N_o$,
non-zero componets of $m$ are $m_{i_1},...,m_{i_n}$ with
$n\in \bf N,$ ${\bar m}:=(m_{i_1},...,m_{i_n})$ for each $m\ne 0$,
$x_{\bar m}:=(x^{i_1},...,x^{i_n})\in {\bf K^n}\hookrightarrow X,$ 
${\bar Q}_0:=1$, where $\bar Q_{\bar m}$ is the Amice basis 
polynomial \cite{ami}.
\par  {\bf 2.1.3.a.} Let $N$ be an analytic manifold modelled on
$Y$ with an atlas 
$$(1)\mbox{ }At(N)=\{ (V_k,\psi _k):\mbox{ }k \in \Lambda _N \} ,
\mbox{ such that }\psi _k: 
V_k\to \psi _k(V_k)\subset Y $$ 
are homeomorphisms, $card(\Lambda _N )\le \aleph _0$ 
and $\theta : M \to N$ be a
$C(t')$-mapping, also $card( \Lambda _M)<\aleph _0$,
where $V_k$ are clopen in $N$, $t'\ge \max (1,t)$ 
is the index of a class of smoothness, that is, for each admissible $(i,j)$:
$$(2)\mbox{ }\theta _{i,j}\in C_*(t',U_{i,j}\to Y)$$
with $*$ empty or an index $*$ taking value $0$ respectively,
$$(3)\mbox{ }\theta _{i,j}:= \psi _i\circ \theta |_{U_{i,j}},$$ where 
$U_{i,j}:=[U_j\cap \theta ^{-1}(V_i)]$ are non-void clopen subsets.
We denote by $C^{\theta }_*(\xi ,M\to N)$ for $\xi =t$ with
$0 \le t \le \infty $ a space of mappings
$f: M\to N$ such that 
$$(4)\mbox{ }f_{i,j}-\theta _{i,j}\in C_*(\xi ,U_{i,j}\to Y).$$
In view of Formulas $(1-4)$ we supply it with an ultrametric 
$$(5)\mbox{ }\rho ^{\xi }_*(f,g)=
\sup_{i,j}\| f_{i,j}-g_{i,j}\| _{C_*(\xi ,U_j\to Y)} $$
for each $0\le \xi < \infty $. 
\par {\bf 2.1.3.b.} Let $M$ and $N$ be two analytic manifolds
with finite atlases, $dim_{\bf K}M=n\in \bf N$, $\theta _{i,j}\in C
(\infty ,U_j\to Y)$ for each $i,j$.
\par We denote by $C_0^{\theta }((t,s),M\to N)$ a completion of a
locally $\bf K$-convex space 
$$(1)\mbox{ } \{ f\in C_0^{\theta }
(t+sn,M\to N):\mbox{ }\rho ^{(t,s)}_0(f,\theta )<\infty $$
$$\mbox{ and for each }\epsilon >0\mbox{ a set }
\{ (k,m): \mbox{ } \sum_{i,j}
|a(m,f^k_{i,j}-\theta ^k_{i,j})| J((t,s),m) >\epsilon  \} 
\mbox{ is finite } \} $$ 
relative to an ultrametric
$$(2)\mbox{ }\rho ^{(t,s)}_0(f,g):=\sup_{i,j,m,k}|a(m,f^k_{i,j}-g^k_{i,j})|
J_j((t,s),m),$$
where $s\in \bf N_o$, $0\le t<\infty $, 
$$(3)\mbox{ }J_j((t,s),m):=\max_{(v\le [t]+ sign(t)+sn)}
\| ({\bar \Phi }^v{\bar Q}_m|_{U_j})(x;$$
$$h_1,...,h_v; \zeta _1,...,\zeta _v)\|_{C_0(0,U_j\to Y)} \mbox{ with}$$ 
$$(4)\mbox{ }h_1=...=h_{\gamma }=e_1,...,h_{(n-1)\gamma +1}=...=
h_{n\gamma }=e_n$$ 
for each integer $\gamma $ such that $1\le \gamma \le s$
and for each $v\in \{ [t]+\gamma n, t+\gamma n \} .$
\par {\bf 2.1.4.} For infinite atlases we use the
traditional procedure of inductive limits of spaces.
For $M$ with the infinite atlas, $card(\Lambda _M)=\aleph _0$, 
and the Banach space $Y$ over $\bf K$ we denote by
$C_*^{\theta }(\xi ,M\to Y)$ for $\xi =t$ with $0\le t\le \infty $
or for $\xi =(t,s)$ a locally $\bf K$-convex space, which is the
strict inductive limit 
$$(1)\mbox{ }C_*^{\theta }(\xi ,M\to Y):=str-ind \{ C_*^{\theta }(\xi ,(U^E
\to Y), \pi ^F_E, \Sigma \} ,$$ 
where $E \in \Sigma $, $\Sigma $
is the family of all finite subsets of $\Lambda _M$ directed by
the inclusion $E<F$ if $E\subset F$, $U^E:=\bigcup_{j\in E}U_j$
(see also \S 2.4 \cite{luseamb}).
\par  For mappings from one manifold into another 
$f: M\to N$ we therefore get the corresponding uniform spaces.
They are denoted by $C_*^{\theta }(\xi ,M\to N).$
\par We introduce notations 
$$(2)\mbox{ }G_i(\xi ,M):=C_0^{\theta }(\xi ,M\to M)\cap Hom(M),$$
$$(3)\mbox{ }Diff(\xi ,M)=C^{\theta }(\xi ,M\to M)\cap Hom(M),$$
that are called groups of 
diffeomorphisms (and homeomorphisms for $0\le t<1$ and $s=0$),
$\theta =id$, $id(x)=x$ for each $x\in M$, 
where $Hom(M):= \{ f:$ $f\in C(0,M\to M),$ $f\mbox{ is bijective },$
$f(M)=M,$ $f\mbox{ and}$ $f^{-1}\in C(0,M\to M) \} $
denotes the usual homeomorphism group.
For $s=0$ we may omit it from the notation, which is always
accomplished for $M$ infinite-dimensional over $\bf K$.
\par {\bf 2.2. Notes.} Henceforth, ultrametrizable separable complete 
manifolds $\bar M$ and $N$ are considered.
Since a large inductive dimension $Ind(\bar M)=0$ 
(see Theorem 7.3.3 \cite{eng}), hence $\bar M$ has not boundaries
in the usual sense. Therefore, 
$$(1)\mbox{ }At({\bar M})=\{(\bar U_j,\bar \phi _j):\mbox{ }
j\in \Lambda _{\bar M} \}$$ has a refinement
$At'(\bar M)$ which is countable and its charts 
$({\bar U'}_j,{\bar \phi '}_j)$ are clopen and
disjoint and homeomorphic with the corresponding balls $B(X,y_j,
{\bar r'}_j)$, where
$$(2)\mbox{ }{\bar \phi '}_j: \mbox{ }{\bar U'}_j\to B(X,{y'}_j,
{\bar r'}_j) \mbox{ for each }j\in {\Lambda '}_{\bar M}$$
are homeomorphisms (see \cite{eng,luum985}).
For $\bar M$ we fix such $At'(\bar M)$.
\par We define topologies of groups $G_i(\xi ,\bar M)$ and locally
$\bf K$-convex spaces $C_*(\xi ,\bar M\to Y)$ relative
to $At'({\bar M})$, where $Y$ is the Banach space over $\bf K$.
Therefore, we suppose also that $\bar M$ and $N$ are
clopen subsets of the Banach spaces $X$ and $Y$ respectively.
Up to the isomorphism of loop semigroups 
(see below their definition) we can suppose that $s_0=0\in \bar M$
and $y_0=0\in N$.
\par For $M={\bar M}\setminus \{ 0\}$ let $At(M)$ consists of charts
$(U_j,\phi _j)$, $j\in \Lambda _M$, while
$At'(M)$ consists of charts $(U'_j,\phi '_j)$, $j\in \Lambda '_M$,
where due to Formulas $(1,2)$ we define
$$(3)\mbox{ }U_1={\bar U}_1\setminus \{ 0\},\mbox{ } 
\phi _1={\bar \phi }_1|_{U_1};\mbox{ }U_j={\bar U}_j\mbox{ and }
\phi _j={\bar \phi }_j\mbox{ for each }j>1,$$ 
$$0\in \bar U_1,\mbox{ }\Lambda _M=\Lambda _{\bar M},\mbox{ }
{U'}_1={\bar U'}_1\setminus \{ 0\},\mbox{ }\phi '_1={\bar \phi '}_1|_{U'_1},
\mbox{ }{U'}_j={\bar U'}_j\mbox{ and }{\phi '}_j={\bar \phi '}_j$$ 
$$\mbox{ for each }j>1,\mbox{ }j\in {\Lambda '}_M=
{\Lambda '}_{\bar M},\mbox{ }{\bar U'}_1\ni 0.$$
\par {\bf 2.3. Definitions and Notes. 1.} Let the spaces 
be the same as in \S 2.1.4 (see Formulas 2.1.4.(1-3))
with the atlas of $M$ defined by Conditions 2.2.(3).
Then we consider their subspaces of mappings preserving marked points:
$$(1)\mbox{ }C^{\theta }_0(\xi ,(M,s_0)\to (N,y_0)):=
\{ f\in C^{\theta }_0(\xi ,\bar M\to N): 
\lim_{|\zeta _1|+...+|\zeta _k|\to 0}{\bar \Phi }^v(f-$$
$$\theta )(s_0;
h_1,...,h_k; \zeta _1,...,\zeta _k)=0\mbox{ for each } 
v\in \{ 0,1,...,[t],t \},\mbox{ } k=[v]+sign \{ v\} \} , $$ 
where for $s>0$ and $\xi =(t,s)$
in addition Condition 2.1.3.b.(4) is satisfied 
for each $1\le \gamma \le s$ and for each
$v\in \{ [t]+n\gamma , t+n\gamma \},$ \\
and the following subgroup:
$$(2)\mbox{ }G_0(\xi ,M):=\{ f\in G_i(\xi ,\bar M): f(s_0)=s_0 \} $$
of the diffeomorphism group, where $s\in \bf N_o$ for 
$dim_{\bf K}M<\aleph _0$ and $s=0$ for $dim_{\bf K}M=\aleph _0$.
\par With the help of them we define the following equvalence relations
$K_{\xi }$: \\
$fK_{\xi }g$ if and only if there exist sequences
$$\{ \psi _n\in G_0(\xi ,M):\mbox{ }n\in {\bf N} \} ,$$
$$ \{ f_n\in C^{\theta }_0(\xi ,M\to N):\mbox{ } n \in {\bf N} \} 
\mbox{ and}$$ 
$$ \{ g_n\in C^{\theta }_0(\xi ,M\to N): n \in {\bf N} \} 
\mbox{ such that}$$ 
$$(3)\mbox{ }f_n(x)=g_n(\psi _n(x))\mbox{ for each }x\in M\mbox{ and }
\lim_{n\to \infty }f_n=f\mbox{ and }\lim_{n\to \infty }g_n=g.$$
Due to Condition $(3)$ these equivalence classes are closed,
since $(g(\psi (x))'=g'(\psi (x))
\psi '(x)$, $\psi (s_0)=s_0$, $g'(s_0)=0$ for $t+s\ge 1$. 
We denote them by $<f>_{K,\xi }$.
Then for $g\in <f>_{K,\xi }$ we write $gK_{\xi }f$ also.
The quotient space $C^{\theta }_0(\xi ,(M,s_0)\to (N,y_0))/K_{\xi }$
we denote by $\Omega _{\xi }(M,N)$, where $\theta (M)=\{ y_0\}$.
\par {\bf 2.3.2.} Let as usually $A\vee B:=A\times \{
b_0\}\cup \{ a_0\}\times B\subset A\times B$ be the wedge 
product of pointed spaces $(A,a_0)$ and $(B,b_0)$, where $A$
and $B$ are topological spaces with marked points $a_0\in A$ and $b_0\in B$.
Then the composition $g\circ f$ of two elements $f, g\in 
C^{\theta }_0(\xi ,(M,s_0)\to (N,y_0))$ is
defined on the domain $\bar M\vee \bar M\setminus \{ s_0\times s_0\}
=:M\vee M$. 
\par Let $M=\bar M\setminus \{ 0\}$ be as in \S 2.2. We fix an 
infinite atlas ${\tilde A}t'(M):=
\{ ({{\tilde U}'}_j,{\phi '}_j): j\in {\bf N} \}$
such that ${\phi '}_j: {{\tilde U}'}_j\to B(X,{y'}_j,{r'}_j)$ 
are homeomorphisms,
$$\lim_{k\to \infty }{r'}_{j(k)}=0 \mbox{ and }
\lim_{k\to \infty }{y'}_{j(k)}=0$$
for an infinte sequence $\{ j(k)\in {\bf N}: k\in {\bf N} \}$ such that
$cl_{\bar M}[\bigcup_{k=1}^{\infty }{{\tilde U}'}_{j(k)}]$ 
is a clopen neighbourhood
of $0$ in $\bar M$, where $cl_{\bar M}A$ denotes the closure of a subset 
$A$ in $\bar M$. In $M\vee M$ we choose the following atlas 
${\tilde A}t'(M\vee M)=\{ (W_l,\xi _l): l\in {\bf N}\}$ 
such that $\xi _l: W_l\to B(X,z_l,a_l)$ are homeomorphisms, 
$$\lim_{k\to \infty }a_{l(k)}=0 \mbox{ and }
\lim_{k\to \infty }z_{l(k)}=0$$ 
for an infinite sequence $\{ l(k)\in {\bf N}: 
k\in {\bf N} \}$ such that $cl_{\bar M\vee \bar M}[
\bigcup_{k=1}^{\infty }W_{l(k)}]$ is a clopen neighbourhood
of $0\times 0$ in $\bar M\vee \bar M$ and 
$$card ({\bf N}\setminus 
\{ l(k): k\in {\bf N} \})=card ({\bf N}\setminus
\{ j(k): k\in{\bf N} \} ) .$$
\par Then we fix a $C(\infty )$-diffeomorphisms $\chi : M\vee M\to M$
such that 
$$(1)\mbox{ }\chi (W_{l(k)})=
{{\tilde U}'}_{j(k)}\mbox{ for each }k\in {\bf N}\mbox{ and}$$ 
$$(2)\mbox{ }\chi (W_l)=
{{\tilde U}'}_{\kappa (l)}\mbox{ for each }l\in ({\bf N}\setminus
\{ l(k): k\in {\bf N} \} ),\mbox{ where}$$
$$(3)\mbox{ }\kappa :({\bf N}\setminus \{ l(k):
k\in {\bf N} \} )\to ({\bf N}\setminus \{ j(k): k\in {\bf N} \} )$$
is a bijective mapping for which
$$(4)\mbox{ }p^{-1}\le a_{l(k)}/{r'}_{j(k)}\le p\mbox{ and }
p^{-1}\le a_l/{r'}_{\kappa (l)}\le p.$$
This induces the continuous injective homomorphism
$$(5)\mbox{ }\chi ^*: C^{\theta }_0(\xi ,
(M\vee M, s_0\times s_0)\to (N,y_0))\to
C^{\theta }_0(\xi ,(M,s_0)\to (N,y_0))\mbox{ such that }$$
$$(6)\mbox{ }\chi ^*(g\vee f)(x)=(g\vee f)(\chi ^{-1}(x))$$
for each $x\in M$, where $(g\vee f) (y)=f(y)$ for $y\in M_2$
and $(g\vee f)(y)=g(y)$ for $y\in M_1$, $M_1\vee M_2=M\vee M$, 
$M_i=M$ for $i=1,2$.
Therefore 
$$(7)\mbox{ }g\circ f:=\chi ^*(g\vee f)$$ 
may be considered as defined on $M$
also, that is, to $g\circ f$ there corresponds the unique element in
$C^{\theta }_0(\xi ,(M,s_0)\to (N,y_0))$.
\par {\bf 2.3.3.} The composition in $\Omega _{\xi }(M,N)$
is defined due to the following inclusion $g\circ f \in 
C^{\theta }_0(\xi ,(M,s_0)\to (N,y_0))$ (see 
Formulas 2.3.2.(1-7)) and then using the equivalence relations
$K_{\xi }$ (see Condition 2.3.1.(3)). 
\par It is shown below that $\Omega _{\xi }(M,N)$
is the monoid, which we call the loop monoid.
\par {\bf 2.4. Note.} For each chart $(V_i, \psi _i)$ 
of $At(N)$ (see Equality 2.1.3.a.(1)) there are local normal coordinates
$y=(y^j: j\in \beta )\in B(Y,a_i,r_i)$, $Y=c_0(\beta ,{\bf K})$.
Moreover, $TV_i=V_i\times Y$, consequently, $TN$ has the disjoint
atlas $At(TN)=$ $\{ (V_i\times X,\psi _i\times I): i \in \Lambda _N\} $,
where $I_Y: Y\to Y$ is the unit mapping, $\Lambda _N\subset \bf N$, $TN$
is the tanget vector bundle over $N$.
\par  Suppose $V$ is an analytic vector field on $N$ (that is,
by definition $V|_{V_i}$ are analytic for each chart and $V\circ
\psi _i^{-1}$ has the natural extension from $\psi
_i(V_i)$ on the balls $B(X,a_i,r_i)$). Then by analogy with the classical
case we can define the following mapping 
$$\bar exp_y(zV)=y+zV(y)\mbox{ for which}$$
$$\partial ^2 \bar exp_y(zV(y))/\partial z^2=0$$ 
(this is the analog
of the geodesic), where $\| V(y)\| _Y|z|\le r_i$ for $y \in V_i$ and
$\psi _i(y)$ is also denoted by $y$, $z\in \bf K$, $V(y)\in Y$.
Moreover, there exists a refinement
$At"(N)=$
$\{ (V"_i,\psi "_i): i\in \Lambda "_N \} $ of $At(N)$. This $At"(N)$ is
embedded into $At(N)$ by charts
such that it is also disjoint and analytic and $\psi "_i(V"_i)$
are $\bf K$-convex in $Y$. 
The latter means that $\lambda x+(1 - \lambda )y\in \psi "_i(V"_i)$
for each $x, y \in \psi "_i(V"_i)$ and each $\lambda \in B({\bf K},0,1)$.
Evidently, we can consider $\bar exp_y$ injective on $V"_i$, $y\in V"_i$.
The atlas $At"(N)$ can be chosen such that 
$$(\bar exp_y|_{V"_i}): V"_i\times B(Y,0,\tilde r_i)\to V"_i$$ 
to be the analytic homeomorphism for each
$i \in \Lambda "_M$, where $\infty >\tilde r_i>0$, $y\in V"_i$,
$$\bar exp_y: (\{ y\} \times B(Y,0,\tilde r_i))\to V"_i$$
is the isomorphism. Therefore, $\bar exp$ is the locally analytic mapping,
$\bar exp: \tilde TN\to N$, 
where $\tilde TN$ is the corresponding neighbourhood
of $N$ in $TN$.
\par  Then 
$$(1)\mbox{ }T_fC^{\theta }_*(\xi ,M\to N)= \{ g
\in C^{(\theta ,0)}_*(\xi ,M\to TN):\mbox{ } \pi _N\circ g=f \} ,$$
consequently, 
$$(2)\mbox{ }C^{\theta }_*(\xi ,M\to TN)=
\bigcup_{f \in C^{\theta }_*(\xi ,M\to N)}T_fC^{\theta }_*(
\xi ,M\to N)=TC^{\theta }_*(\xi ,M\to N),$$ 
where $\pi _N: TN \to N$ is the natural projection,
$*=0$ or $*=\emptyset $ ($\emptyset $ is omitted).
Therefore, the following mapping 
$$(3)\mbox{ }\omega _{\bar exp}: \mbox{ }
T_fC^{\theta }_*(\xi ,M\to N) \to 
C^{\theta }_*(\xi ,M\to N)$$ is defined 
by the formula given below
$$(4)\mbox{ }\omega _{\bar exp}(g(x))=
{\bar exp}_{f(x)}\circ g(x),$$ 
that gives charts
on $C^{\theta }_*(\xi ,M\to N)$ induced by charts on 
$C^{\theta }_*(\xi ,M\to TN)$.
\par {\bf 2.5. Definition and Note.} In view of Equalities 2.4.(1,2)
the space $C^{\theta }_0(\xi ,\bar M\to N)$ is isomorphic with
$C^{\theta }_0(\xi ,(M,s_0)\to (N,y_0))\times N^{\xi }$, where $y_0=0$
is the marked point of $N$. Here 
$$(1)\mbox{ }N^{\xi }:=N\otimes (\bigotimes_{j=1}^d
\tilde L_{\xi }(X^j\to Y))\mbox{ for }t\in {\bf N_0}
\mbox{ with }t+s>0;$$
$$(2)\mbox{ }N^{\xi }=N\mbox{ for }t+s=0;$$ 
$$(3)\mbox{ }N^{\xi }=N\otimes (
\bigotimes_{j=1}^d\tilde L_{\xi }(X^j\to Y))\otimes C^0_0(0,M^k\to 
Y_{\bf \lambda })\mbox{ for }t\in {\bf R}\setminus {\bf N},$$
where $N^{\xi }$ is with the product topology,
$d=[t]$ for $\xi =t$, $d=[t]+n\alpha $ 
for $\xi =(t,s)$ with $\alpha =dim_{\bf K}M<\aleph _0$, 
when $s>0,$ $k=d+sign \{ t\}$,
$Y_{\bf \lambda }:=c_0(\beta ,{\bf \lambda })$, ${\bf \lambda }$
is the least subfield of $\bf \Lambda _p$ such that
${\bf \lambda }\supset {\bf K}\cup j_{ \{ t \} } ({\bf K})$ (see
\S 2.1 \cite{lubp2}).
Then $\tilde L_{\xi }(X^j\to Y)$ denotes the Banach space of 
continuous $j$-linear operators $f_j: X^j\to Y$ with 
$$(4)\mbox{ }\| f_j\|_{\tilde L_{\xi }
(X^j\to Y)}:=\sup_{i,m}\| f^i_j\|_m \mbox{ and }$$
$$(5)\mbox{ }\lim_{i+|m|+k
\to \infty }\| f^i_j\|_m=0,\mbox{ where}$$ 
$$(6)\mbox{ }\| f^i_j\|_m:=
\sup_{0\ne h_l\in {\bf K}^k, l=1,...,j} \| f^i_j(h_1,...,h_j) \|_Y
J'(\xi ,m)/(\| h_1\|_X...\| h_j\|_X),$$ 
${\bf K}^k:=sp_{\bf K}(e_1,...,e_k)\hookrightarrow X$
is a $\bf K$-linear span of the standard basic vectors, 
$m=(m_1,...,m_k)$, $|m|=m_1+...+m_k$, $k\in \bf N$;
$h_1=...=h_{m_1},$...,$h_{m_{k-1}+1}=...=h_{m_k}$ for $s=0$;
in addition Condition 2.1.3.b.(4) is satisfied for each
$0<\gamma \le s$, when $s>0$; $f=(f_0,f_1,...,f_j,...)\in N^{\xi }$,
$\sum_if^i_jq_i=f_j$, 
$f^i_j: X^j\to \bf K$, 
$$J'(\xi ,m):=|\partial ^m{\bar Q}_m(x)|_{x=0}|_{\bf K}$$
(see \S 2.2 \cite{lubp2} and Equations 2.1.2.(1-5), 2.1.3.b.(1-3)).
\par {\bf 2.6. Definitions.} A function $f: {\bf K}\to \bf C$ is called
pseudo-differentiable of order $b$, if there exists the following integral:
$$(1)\mbox{ }PD(b,f(x)):=\int_{\bf K} [(f(x)-f(y)) \times g(x,y,b)]
v(dy),$$ 
where $g(x,y,b):=\mid x-y\mid ^{-1-b}$
with the nonnegative Haar measure $v$ and
$b \in {\bf C}$ (see also \S 2.1 \cite{lubp2}).
We introduce the following notation $PD_c(b,f(x))$ for
such integral by $B({\bf K},0,1)$ instead of the entire $\bf K$.
\par {\bf 2.7. Definitions.} Let $G$ be a topological Hausdorff semigroup
and $({\sf M},{\bf R})$ be a space $\sf M$ 
of measures on $(G,Bf(G))$ with values
in $\bf R$, where $Bf(G)$ denotes the Borel $\sigma $-algebra of 
$G$. Let also 
$G'$ and $G"$ be dense subsemigroups in $G$
such that $G"\subset G'$ and a topology $\sf T$ on
$\sf M$ is compatible with $G'$, that is, $\mu \mapsto \mu _h$
is the homomorphism
of $({\sf M},{\bf R})$ into itself for each $h \in G'$,
where $\mu _h(A):=\mu (h\circ A)$ for each $A\in Bf(G)$. 
Let $\sf T$  be the topology of convergence for each $E \in Bf(G)$.
If $\mu \in ({\sf M},{\bf R})$ and $\mu _h\sim \mu $ 
are the equivalent measures
for each $h\in G'$ then $\mu $ is called quasi-invariant on 
$G$ relative to $G'$. We shall consider $\mu $ with the continuous 
quasi-invariance factor 
$$(1)\mbox{ }\rho _{\mu }(h,g):=\mu _h(dg)/\mu (dg).$$
If $G$ is a group, then we use the traditional definition of
$\mu _h$ such that $\mu _h(A):=\mu (h^{-1}\circ A)$.
\par Let $S(r,f)=g(r,f)$ be a curve on the
subsemigroup $G"$, such that $S(0,f)=f$ and there exists
$\partial S(r,f)/\partial r\in TG"$
and $\partial S(r,f)/\partial r|_{r=0}=:A_f\in T_fG"$,  where
$r\in B({\bf K},0,R)$, $\infty >R\ge 1$.
Then a measure $\mu $ on $G$ is called pseudo-differentiable of order $b$ 
relative to $S$ if there exists $PD_c(b,\bar S(r,\mu )(B))$ by $r\in
B({\bf K},0,1)$ for each $B \in Bf(G)$,
where $\bar S(r,\mu )(B):=
\mu (S(-r,B))$ for each $B \in Bf(G)$. A measure
$\mu $ is called pseudo-differentiable of order $b$ if there exists a dense
subsemigroup $G"$ of $G$ such that $\mu $ is pseudo-differentiable
of order $b$ for each curve $S(r,f)$ on $G"$ described above,
where $b\in \bf C$.
\par Naturally Definitions 2.7 have generalizations, when $G$ is a 
topological manifold on which a topological group 
(or a semigroup) $G'$ acts continuously
from the left $G'\times G\ni (g,x)\mapsto gx\in G.$
\par {\bf 2.8.} {\bf Note.} Now let us describe dense
loop submonoids which are necessary for the investigation
of quasi-invariant measures on the entire monoid. 
For finite $At(M)$ and $\xi =(t,s)$
let $C_{0, \{ k\} }^{\theta }(
\xi ,M\to Y)$ be a subspace of $C_0^{\theta }(\xi ,M\to Y)$
consisting of mappings $f$ for which 
$$(1)\mbox{ }\| f- \theta\|_{
C_{0, \{ k\} }^{\theta }( 
\xi ,M\to Y)}:=\sup_{i,m,j}|a(m,f^i|_{U_j})|_{\bf K}J_j(\xi ,m)p^{k(i,m)}
<\infty \mbox{ and}$$ 
$$(2)\lim_{i+|m|+Ord(m)\to \infty } \sup_j
|a(m,f^i|_{U_j})|_{\bf K}J_j(\xi ,m)p^{k(i,m)}=0,$$ 
where $k(i,m):=c'\times i+c\times (|m|+Ord(m))$, $c'$ and $c$ are
non-negative constants, $|m|:=\sum_im_i$,
$$Ord(m):=\max \{ i:\mbox{ }m_i>0\mbox{ and }m_l=0\mbox{ for each }
l>i \} $$
(see also Formulas 2.1.2.(2) and 2.1.3.b.(3)).
\par For finite-dimensional $M$ over $\bf K$
this space is 
isomorphic with $C_{0, \{ k'\} }^{\theta }( \xi ,M\to Y)$,
where $k'(i,m)=c'\times i+c\times |m|$. 
For finite-dimensional $Y$ over $\bf K$ the space
$C_{0, \{ k\} }^{\theta }( \xi ,M\to Y)$  is isomorphic with
$C_{0, \{ k"\} }^{\theta }( \xi ,M\to Y)$, where $k"(i,m)=c\times
(|m|+Ord(m))$. For $c'=c=0$ this space coincides with
$C_0^{\theta }(\xi ,M\to Y)$ and we omit $\{ k\}$.
\par Then as in \S 2.3 we define spaces
$C_{0, \{ k\} }^{\theta }( \xi ,(M,s_0)\to 
(N,0))$, groups
$$(3)\mbox{ }G^{ \{ k\} }(\xi ,M):=C^{id}_{0, \{ k\} }( 
\xi , M\to M)\cap Hom(M),$$ 
$$(4)\mbox{ }G^{ \{ k\} }_0(\xi ,M):=
\{ \psi \in G^{ \{ k\} }(\xi ,M): \psi (s_0)=s_0 \} $$
and the equivalence relation $K_{\xi , \{ k\} }$
in it for each $M$ and $N$ from \S 2.1 and \S 2.2. 
Therefore, 
$$(5)\mbox{ }G':=\Omega _{\xi }^{ \{ k\} }
(M,N)=:C_{0,\{ k\} }^0(\xi ,(M,s_0)\to (N,0))/K_{\xi ,\{ k\} } $$
is the dense submonoid in $\Omega _{\xi }(M,N)$.
\par {\bf 2.9. Note and Definition.} For a commutative 
monoid $\Omega _{\xi }(M,N)$
with the unity and the cancellation property (see \cite{lubp2})
there exists a commutative group $L_{\xi }(M,N)$
equal to the Grothendieck group. This group is the quotient group
$F/\sf B$, where $F$ is a free Abelian group generated by 
$\Omega _{\xi }(M,N)$ and $\sf B$ is a closed subgroup of $F$ generated by
elements $[f+g]-[f]-[g]$, $f$ and $g\in \Omega _{\xi }(M,N)$,
$[f]$ denotes an element of $F$ corresponding to $f$. The natural mapping 
$$(1)\mbox{ }\gamma : \Omega _{\xi }(M,N)\to L_{\xi }(M,N)$$ 
is injective.
We supply $F$ with a topology inherited from
the Tychonoff product topology of $\Omega _{\xi }(M,N)^{\bf Z}$,
where each element $z$ of $F$ is 
$$(2)\mbox{ }z=\sum_fn_{f,z}[f],$$
$n_{f,z}\in \bf Z$ for each $f\in \Omega _{\xi }(M,N)$,
$$(3)\mbox{ }\sum_f|n_{f,z}|<\infty .$$ 
In particular $[nf]-n[f]\in \sf B$, where 
$1f=f$, $nf=f\circ (n-1)f$ for each 
$1<n\in \bf N$, $f+g:=f\circ g$. 
We call $L_{\xi }(M,N)$ the loop group. 
\par {\bf 2.10. Note.}
Let $\Omega _{\xi }^{\{ k\} }(M,N)$ be the loop submonoid as in \S 2.8
such that $c>0$ and $c'>0$.
Then it generates the loop group $G':=L _{\xi }^{\{ k\} }(M,N)$
as in \S 2.9 such that $G'$ is the dense subgroup in $G=L _{\xi }(M,N)$.
\par {\bf 2.11. Remarks.} Let $M$ be a manifold on the Banach space $X$ 
with an atlas $At(M)$ consisting of disjunctive charts 
$(U_j,\phi _j)$, $j\in \Lambda $,
$\Lambda \subset \bf N$, where $U_j$ and $\phi _j(U_j)$,
are clopen in $M$ and $X$ respectively, $\phi _j:
U_j\to \phi _j(U_j)$ is a homeomorphism, 
also $\phi _j(U_j)=B(X,x_j,r_j)$ is a ball in $X$ with
a radius $0<r_j<\infty $ for each $j$.
\par For $\Lambda =\omega _0$ we define a Banach space 
$${\tilde C}_*(t,M\to X):= \{ f|_{U_j}\in C_*(t,U_j\to X),  
\| f\|_{C_*(t,M\to X)}:=\sup_{j\in \Lambda }
(\| f|_{U_j}\|_{C_*(t,U_j\to X)}$$
$$/\min (1,r_j))<\infty \mbox{ and }
(\| f|_{U_j}\|_{C_*(t,U_j\to X)})
/\min (1,r_j))\to 0 \mbox{ while }j\to \infty \} ,$$ 
where $0\le t< \infty $,
$*=0$ for spaces $C_0(t,U\to X)$, $*=\emptyset $
or simply is omitted for $C(t,U\to X)$. 
For the finite atlas $At(M)$ the spaces
${\tilde C}_*(t,U\to X)$ and $C_*(t,U\to X)$ are 
linearly topologically isomorphic.
By $C_*^{\theta }(t,M\to M)$ for $0\le t\le \infty $
is denoted the following space of functions
$f: M\to M$ such that $(f_i-\theta _i)\in
C_*(t,M\to X)$ for each $i\in \Lambda $ and $f_i=\psi _i
\circ f$, $\theta _i=\psi _i\circ \theta $.
We introduce the following group 
$$G(t,M):={\tilde C}_0^{id}(t,M\to M)\cap Hom(M),$$
which is called the diffeomorphism group 
(and the homeomorphism group for $0\le t<1$), where $Hom(M)$
is the group of continuous homeomorphisms.
\par Each function $f\in C_0(t,M\to X)$ 
has the following decomposition:
$$f(x)|_{U_j}=\sum_{(i\in {\bf N}, n \in {\bf N_o} )}
f^i(n;x)|_{U_j}e_i\tilde z(n),\mbox{ and  }
\{ e_i\tilde z(n)(\bar Q_m(x)|_{U_j}): $$
$i,n, Ord(m)=n, j \}$ is the orthogonal basis, moreover, 
$$f_n(x)|_{U_j}:=\sum_if^i(n;x)|_{U_j}e_i \in 
C_0(t,U_j\to X),\mbox{ where}$$ 
$$X_{\tilde z(n)}:=\{ f_n(x): f_n|_{U_j}\in C_0(t,U_j\to X) \} $$
is the Banach space with the norm 
induced from $C_0(t,M\to X)$ such that
$$f^i(n;x)|_{U_j}:=\sum_{(Ord \mbox{ } m=n,
m=(m(1),...,m(n)),m(j) \in {\bf N_o} )} a(m,f^i|_{U_j})
\bar Q_m(x)|_{U_j},$$ where
$\bar Q_m(x)|_{U_j}=0$ for $x\in M\setminus U_j$.
\par For the manifold $M$ we fix a subsequence
$\{ M_n: n\in {\bf N_o} \} $ of submanifolds
in $M$ such that $M_n\hookrightarrow M_{n+1}
\hookrightarrow ... M$ for each $n$, $dim_{\bf K}M_n=\beta (n)\in \bf N$
for each $n\in \bf N_o$, $\bigcup_nM_n$ is dense in $M$, where
$\beta (n)<\beta (n+1)$ for each $n$ and there exists
$n_0\in \bf N$ with $\beta (n)=n$ for each $n>n_0$.
\par We take the following subgroup
$$G':=\{f \in G(t,M):
(f^i(n;x)-id^i(n;x))=:g^i(
n;x) \in C_0(t_n,M_n\to {\bf K}) 
\mbox{ and}$$
$$\mid a(m;g^i(n;x)|_{U_j})
\mid J_j(t_n,m)\le c(f)p^{v'(m,j,i)} \} ,$$
where $c(f)>0$ is a constant, $v'(m,j,i)=-c'i-c'n-c"j $, $n=Ord (m)$,
$c'=const >0$ and $c"=const \ge 0$, $c">0$ for $\Lambda =\omega _0$,
$t_n=t+s(n)$ for $0\le t<\infty $, $s(n)> n$ for each $n$
and $\liminf_{n\to \infty }s(n)/n=: \zeta >1$.
Then there exists the following ultrametric in $G'$:
$$d(f,id)
=\sup _{m,n,j}\{ \mid a(m;g^i(n;x)|_{U_j})\mid J_j(t_n,m)
p^{-v'(m,j,i)} \} $$. 
\par {\bf 2.12. Note.} At first it is necessary to prove theorems about the
quasi-invariance and the pseudo-differentiability of transition measures
of stochastic processes on Banach spaces over local fileds.
We consider two types of measures on $c_0(\omega _0,{\bf K}).$
The first is the $q$-Gaussian measure 
$$\mu =\mu _{J,\gamma,q}
:=\bigotimes_{j=1}^{\infty }\mu _j(dx^j),\mbox{ where }
\mu _j(dx^j)=C_{|\zeta _j|^{-q},\gamma _j,q}
f_{|\zeta _j|^{-q},\gamma _j,q}v(dx^j)$$
(see \S 2 \cite{lunast2}).
The characteristic functional of the $q$-Gaussian measure 
is positive definite, hence $\mu $ is nonnegative
(see also \S 2.6 \cite{lu6}).
The second is specified below and is the particular case
of measures considered in \S 4.3 \cite{lunast1}.
\par Let $w$ be the real-valued nonnegative Haar measure on $\bf K$
with $w(B({\bf K},0,1))=1$.
We consider the following measure $\mu $ on $c_0(\omega _0,{\bf K})$  \\
$(i)\quad \mu (dx)=\bigotimes_{j=1}^{\infty }\mu _j(dx^j),$ where 
$x\in c_0(\omega _0,{\bf K})$, $x=(x^j: j\in \omega _0),$ 
$x^j\in \bf K$, $x=\sum_jx^je_j$, $e_j$ is the standard othonormal base in 
$c_0(\omega _0,{\bf K})$.
\par Let now on the Banach space $c_0:=c_0(\omega _0,{\bf K})$ 
there is given an operator $J\in L_1(c_0)$ such that
$Je_i=v_ie_i$ with $v_i\ne 0$ for each $i$. We consider a measure 
$\nu _i(dx):=f_i(x)w(dx)$ on $\bf K$,
where $f_i: {\bf K}\to [0,1]$ is a function belonging to the space
$L^1({\bf K},w,{\bf R})$ such that 
$f_i(x)=f(x/v_i)+h_i(x/v_i)$, 
where $f$ is a locally constant positive function,
$f(x)=\sum_{j=1}^{\infty }C_jCh_{B_j}(x)$, 
$B_j:=B({\bf K},x_j,r_j)$ is a ball in $\bf K$, 
$Ch_V$ is the characteristic function of a subset  $V$ in $\bf K$, 
that is, $Ch_V(x)=1$ for each $x\in V$, $Ch_V(x)=0$ for each 
$x\in {\bf K}\setminus V$, $x_1:=0$, $r_1:=1$,
$\inf_jr_j=1$, $\{ B_j: j \} $ is the disjoint covering of $\bf K$,
$1\ge C_j>0$, $\lim_{|x|\to \infty }f(x)=0$, $h_i\in 
L^1({\bf K},w,{\bf R})$ such that $ess_w-\sup_{x\in {\bf K}}
|h_i(x)/f(x)|=\delta _i < 1$, $\sum_i\delta _i<\infty $
and $\nu _i({\bf K})=1$. 
Then $\nu _i(S)>0$ for each open subset $S$ in $\bf K$.
There exists a $\sigma $-additive product measure \\
$(ii)\quad \mu _J(dx):=\prod_{i=1}^{\infty }\mu _i(dx^i)$ 
on the $\sigma $-algebra of Borel subsets of $c_0$, 
since the Borel $\sigma $-algebras defined for the
weak topology of $c_0$ and for the norm topology of $c_0$ coincide,
where $\mu _i(dx^i):=\nu (dx^i/v_i)$.
\par Let $A: c_0\to c_0$ be a linear topological isomorphism,
that is, $A$ and $A^{-1}\in L(c_0)$, then 
for a measure $\mu $ on $c_0$ there exists its image
$\mu _A(S):=\mu (A^{-1}S)$ for each Borel subset $S$ in $c_0$.
In view of Proposition 2.12.2 \cite{lunast1} $L_q(c_0)$ is 
the ideal in $L(c_0)$. This produces new $q$-Gaussian measures 
$(\mu _{J,\gamma ,q})_A=:\mu _{AJ,A^*\gamma ,q}$
and measures of the second type $(\mu _J)_A=:\mu _{AJ}$.
In view of \S 2.9 \cite{lunast1} each injective 
linear operator $S\in L_q(c_0)$ with $E(c_0)$ dense in $c_0$
can be presented in the form $S=AJ$. Hence for each
such $S$ there exists the $\sigma $-additive measure
$\mu _{S,S^*\gamma ,q}$ and $\mu _S$. These measures 
are induced by the corresponding
cylinder measures $\mu _{I,\gamma ,q}$ or $\mu _I$
on ${\bf K}^{\aleph _0}$, where $I$ is the unit operator, since
$c_0$ in the weak topology is isomorphic with ${\bf K}^{\aleph _0}$.
Here the algebra $\sf U$ of cylindrical subsets is generated by
subsets $\pi _V ^{-1}(A)$, where $A$ is a Borel subset in $\bf K^n$,
$card (V)=n <\aleph _0$, $V\subset \bf N$, $\pi _V: 
{\bf K}^{\aleph _0}\to \prod_{i\in V}{\bf K}_i$ is the natural projection.
\par On the space $C^0_0(T,H)=C^0_0(T,{\bf K})\otimes H$ 
let $S=S_1\otimes S_2$ and $\gamma =\gamma ^1\otimes \gamma ^2$, 
where $S_1$ is a linear operator on
$C^0_0(T,{\bf K})$ and $S_2$ is a linear operator on $H$,
$\gamma ^1\in C^0_0(T,{\bf K})$, $\gamma ^2\in H$
such that the measure $\mu _{S,\gamma ,q}$ is the product of
measures $\mu _{S_1,\gamma ^1,q}$ on $C^0_0(T,{\bf K})$ and
$\mu  _{S_2,\gamma ^2,q}$ on $H$, analogously $\mu _S$ is the product 
of measures $\mu _{S_1}$ on $C^0_0(T,{\bf K})$ and $\mu _{S_2}$ on $H$.
With the help of such measures on the space
$C^0_0(T,H)$ the stochastic process $w(t,\omega )$ is defined 
as in \S \S 4.2 and 4.3 \cite{lunast1} and \S 3.2 \cite{lunast2}.
\par {\bf 2.13.} Let $Y$ be a Banach space over the local field
$\bf K$ and $V$ be a neighbourhood of zero in $Y$. 
Consider either the measure $\mu _{S,\gamma ,q}$ or $\mu _S$
outlined in \S 2.12. Suppose that
in stochastic antiderivational equations 3.4.(i) and 3.5.(i) \cite{lunast2}
mappings $a$ and $E$ be dependent on the parameter $y\in V$, that is,
$a=a(t,\omega ,\xi ,y)$ and $E=E(t,\omega ,\xi ,y)$;
moreover, $a_{k,l}=a_{k,l}(t,\xi ,y)$ for each $k$ and $l$ in the 
latter equation, the condition $3.4.(LLC)$ \cite{lunast2}
is satisfied for each $0<r<\infty $
with the constant $K_r$ independent from $y\in V$
for each $y\in V$. 
Evidently, Equation 3.4.(i) is the particular case of 3.5.(i),
when in the latter equation the corresponding $a_{0,1}$ and
$a_{1,0}$ are chosen with all others $a_{k,l}=0$ (when $k+l\ne 1$).
Let also 
\par $(i)$ $a,$ $E$ and $a_{k,l}$ be of class $C^1$ by $y\in V$
such that \\
$a\in C^1(V,L^q(\Omega ,{\sf F},\lambda ;C^0(B_R,L^q(
\Omega ,{\sf F},\lambda ;C^0(B_R,H)))))$ and  \\
$E\in C^1(V,L^q(\Omega ,{\sf F},\lambda ;C^0(B_R,L(L^q(
\Omega ,{\sf F},\lambda ;C^0(B_R,H)))))),$ \\
$a_{m-l,l} \in C^1(V,
C^0(B_{R_1}\times B(L^q(\Omega ,{\sf F},\lambda ;C^0(B_R,H)),
0,R_2),L_m(H^{\otimes m};H)))$ (continuous and bounded on its domain)
for each $n, l,$ $0<R_2<\infty $ and \\
$\lim_{n\to \infty } \sup_{0\le l\le n}\|a_{n-l,l}
\|_{C^1(V,C^0(B_{R_1}\times B(L^q(\Omega ,{\sf F},\lambda ;C^0(B_R,H)),
0,R_2),L_n(H^{\otimes n},H)))}=0$ 
for each $0<R_1\le R$ when $0<R<\infty $, or each $0<R_1<R$
when $R=\infty $, for each $0<R_2<\infty $;
\par $(ii)$ $ker (E(t,\omega ,\xi ,y))=0$ for
each $t$, $\xi $ and $y$, also for $\lambda $-almost every $\omega $;
\par $(iii)$ $a_y(t,\omega ,\xi ,y)$ 
and $\partial a(t,\omega ,\xi ,y)/\partial y\in
X_{0,s}(H):=\{z: S^{-1}z\in H_s \} $ 
and $\partial E(t,\omega ,\xi ,y)/\partial y
\in L_r(H)$ for $\lambda $-almost all 
$\omega $ and each  $t$, $\xi $, $y$, where $H_s:= \{ z : z\in H; 
\sum_{j=1}^{\infty } |z_j|^s <\infty \} $ 
for each $0<s<\infty $, $H_{\infty }:=H$, with 
$s=r=q$ for $\mu _{S,\gamma ,q}$; $s=\infty $ and $r=0$ for
the measure of the second type $\mu _S$, $z_j$ 
are the coordinates of the vector $z$ in the standard base
in $H$; in addition for Equation 3.5.(i)
\par $(iv)$ $\partial a_{l,k}(t,\omega ,\xi ,y)/
\partial y\in L_{k+l,r}(H^{\otimes (k+l)};H)$ 
for each $l$ and each $k$ with either 
$r=q$ or $r=0$ correspondingly.
The following theorem states the quasi-invariance of 
the transition measure
$\mu ^{F_{t,t_0}} (\{ \omega : 
\xi (t_0,\omega ,y)=0, \xi (t,\omega ,y)\in A \} )=:
P_y(A),$ where $F_{t,u}(\xi ):=\xi (t,\omega ,y)- \xi (u,\omega ,y)$.
\par {\bf Theorem.} {\it Let either Conditions $(i-iii)$ 
or $(i-iv)$ be satisfied, then the transition measure 
$P_y(A)$ 
of the stochastic process $\xi (t,\omega ,y)$ being the solution 
of Equation either 3.4(i) or 3.5.(i) \cite{lunast2}
and depending on the parameter $y\in V$ is quasi-invariant relative to 
each mapping $U(y_2,y;\xi (t,\omega ,y)):=\xi (t,\omega ,y_2)$
for each $y$ and $y_2\in V$.}
\par {\bf Proof.} The Kakutani theorem  (see II.4.1 \cite{dalf})
states, whether $\prod_{k=1}^{\infty }\alpha _k$ converges 
to a positive number or diverges to zero, 
the measure $\mu $ is absolutely continuous or orthogonal
with respect to $\nu $, correspondingly, where $\alpha _k:=
\int_{X_k}(p_k(x_k))^{1/2}\nu _k(dx_k)$, $\mu _k$ is absolutely 
continuous relative to $\nu _k$,
$\mu =\otimes _k\mu  _k$, $\nu =\otimes _k\nu _k$, $\mu _k$ and $\nu _k$ 
are probability measures on measurable spaces $X_k$ for each $k\in \bf N$,
$p_k(x):=\mu _k(dx)/\nu _k(dx)$. In the first case $\prod_kp_k(x_k)$
converges in the mean to $\mu (dx)/\nu (dx)$. 
In the considered here case let $X_k=\bf K$ for each $k\in \bf N$.
Let $\mu _k(dx)=Cf(x-y)v(dx)$, $\nu _k(dx)=Cf(x)v(dx)$, where
$v$ is the non-negative Haar measure on $\bf K$, $f$ is a 
positive function such that $f\in L^1({\bf K},v,{\bf R})$,
$C=const>0$ such that $\nu ({\bf K})=1$.
Then $p_k(x)=f(x-y)/f(x)$ and $\alpha _k=\int_{\bf K}
(f(x-y)f(x))^{1/2}v(dx).$ For the $q$-Gaussian measure
$f(x)=\int_{\bf K}exp(-\beta |x|^q)\chi _{\gamma }(x)\chi _1(-zx)v(dx)$
(see \cite{roo} and \S 7 \cite{vla3}).
If $|yx|\le 1$, then $\chi _1(yx)=1$. Therefore,
there is a constant $C_1>0$ independent from $\beta $ and
$\gamma $ such that
$|f(z-y)-f(z)|\le |f(z)|(1+C_1exp(-\beta r^{-q}))$
for each $y$ with $|y|<r$,
where $\beta r^{-q}>1$, since due to Cauchy-Schwarz-Bunyakovskii inequality
$$|\int_{|x|>1/r}exp(-\beta |x|^q)\chi _{\gamma }(x)\chi _1
(-(z-y)x)v(dx)|\le $$ 
$$|\int_{|x|>1/r}exp(-\beta |x|^q)\chi _{\gamma }
(x)\chi _1(-zx)v(dx)| g(y,z) \le |f(z)|g(y,z),$$ 
where $g(y,z):=|\int_{|x|>1/r}exp(-\beta |x|^q)
\chi _{\gamma }(x)\chi _1(-zx)\chi _1(2yx)v(dx)|$.
Let $|y_j/v_j|=:r_j<1$ for each $j>j_0$, 
then $|\alpha _j-1|\le Cexp(-\beta _jr_j^{-q})$
for each $j>j_0$, where $C=const >0$.
In view of Proposition 2.12.2 \cite{lunast1} and the Kakutani theorem
$\mu ^z_{S,\gamma ,q}$ is equvalent to $\mu _{S,\gamma ,q}$
for each $z\in X_{0,q}(C^0_0(T,H))$, where
$\mu ^z(A):=\mu (A-z)$ for each Borel subset $A$ in $C^0_0(T,H)$, 
that is, $\mu _{S,\gamma ,q}$ is quasi-invariant 
relative to shifts $z\in X_{0,q}(C^0_0(T,H))$.
\par For the measure $\mu _J$ and $|y|<1/|v|$
there is the equality $f((x-y)/v)=f(x/v)$
for each $x\in \bf K$ and $0\ne v\in \bf K$.
In view of the definition of $f_k$ there is the equality
$p_k(x)=f_k(x-y_k)/f_k(x)=[f((x-y_k)/v_k)/f(x/v_k)]
[1+h_k((x-y_k)/v_k)/f((x-y_k)/v_k)]/ [1+h_k(x/v_k)/f(x/v_k)]$.
If $|y_k/v_k|\le 1$, then $f((x-y_k)/v_k)/f(x/v_k)=1$
for each $x\in \bf K$. From the conditions imposed on
$h_k$ and $f$ and the Kakutani theorem  and Proposition
2.12.2 \cite{lunast1} it follows, that $\mu _S$
is quasi-invariant relative to shifts $z\in X_{0,\infty }(C^0_0(T,H))$.
\par The  quasi-invariance factor $\rho (z,x):=\mu ^z(dx)/\mu (dx)$ 
is Borel measurable  as  follows from the construction  of $\mu $
and the Kakutani theorem and the Lebesgue theorem about
majorized convergence (see \S 2.4.9 \cite{feder}),  
since this is true for each its one-dimensional projection.
The Banach theorem states: if $G$ is a topological group and
$A\subset G$ is a Borel measurable set of second  category,
then $A\circ A^{-1}$ is a neighbourhood of unit
(see \S 5.5 \cite{chris}).
The quasi-invariance factor satisfies the cocycle condition:
$\rho (z+h,x)=\rho (z,x-h)\rho (h,x)$ for each $z$ and $h\in X_{0,s}
(C^0_0(T,H))$ and each $x\in C^0_0(T,H)$.
Therefore, in view of the Lusin theorem (see \S 2.3.5 \cite{feder})
$\rho (z,x):=\mu ^z(dx)/\mu (dx)$ is such that
$\mu (W_L)=1$ for each finite-dimensional subspace
$L$ in $X_{0,s}(C^0_0(T,H))$, 
where either $\mu =\mu _{S,\gamma ,q}$ or $\mu =\mu _S$, 
$W_L:=\{ x: \rho (z,x) $ is defined and continous by $z\in  L \} $.
\par In view of the preceding consideration $\lim_{n\to \infty }
\rho (\hat P_nz,x)=\rho (z,x)$ for $\mu $-almost all $x\in C^0_0(T,H)$,
moreover, this convergence is uniform by $z$ in each ball
$B(L,0,c)$ for each finite-dimensional subspace $L$ in $X_{0,s}
(C^0_0(T,H))$,
where $\hat P_n$ is a projection on a subspace $sp_{\bf K}(e_1,...,e_n)
=\bf K^n$, where $\{ e_j: j\} $ is the orthonormal base
in $X_{0,s}(C^0_0(T,H))$. Evidently, 
$X_{0,s}(C^0_0(T,H))$ is dense in $C^0_0(T,H).$
\par Stochastic antiderivational Equation 3.4.(i) \cite{lunast2}
is the particular case of 3.5.(i). Therefore, it is sufficient 
to consider the latter equation. Below it is shown, that
the one-parameter family of solutions $\xi (t,\omega ,y)$
is of class $C^1$ by $y\in V$.
Let $X_0(t,y)=x(y)$,..., 
\par $$X_n(t,y)=x(y)+
\sum_{m+b=1}^{\infty }\sum_{l=0}^m(
{\hat P}_{u^{b+m-l},w(u,\omega )^l}[a_{m-l+b,l}
(u,X_{n-1}(u,\omega ,y),y)\circ $$
$$(I^{\otimes b}\otimes a^{\otimes (m-l)}
\otimes E^{\otimes l})])
|_{u=t},$$ consequently,
$$X_{n+1}(t,y)-X_n(t,y)=$$
$$\sum_{m+b=1}^{\infty }\sum_{l=0}^m(
{\hat P}_{u^{b+m-l},w(u,\omega )^l}[a_{m-l+b,l}
(u,X_n(u,y),y)-a_{m-l+b,l}(u,X_{n-1}(u,y),y)]$$
$\circ (I^{\otimes b}\otimes a^{\otimes (m-l)}\otimes E^{\otimes l})])
|_{u=t},$ \\
where $t_j=\sigma _j(t)$ for each $j=0,1,2,...$,
for the shortening of the notation $X_n,$ $x$ and $a_{l,k}$ are written 
without the argument $\omega $, $a$ and $E$ are written 
without their variables. Then
$$M \sup_y \| {\hat P}_{u^{b+m-l},w(u,\omega )^l}[a_{m-l+b,l}
(u,X_n(u,y),y)-$$ $$a_{m-l+b,l}(u,X_{n-1}(u,y),y)
]|_{(B_{R_1}\times B(L^q,0,R_2)\times V)}\circ ( $$ 
$$I^{\otimes b}\otimes a^{\otimes (m-l)}\otimes E^{\otimes l})])
|_{u=t} \| ^g \le K(M \| {\hat P}_{u^{b+m-l},w(u,\omega )^l} \| ^g)
\| a_{m-l+b,l} |_{(B_{R_1}\times B(L^q,0,R_2)\times V)}\| ^g$$ 
$$ (M \sup_{u,y} \| X_n(u,y)-X_{n-1}(u,y) \| ^g)
(M\sup_{u,y} \| a \| ^{m-l})(M\sup_{u,y} \| E \| ^l),$$
where $X_n\in C^0_0(B_R,H)$ for each $\omega $, $y\in V$ and
for each $n$, $K$ is the same constant as in \S 3.4, $1\le g<\infty $. 
On the other hand,
$$X_1(t,y)=x(t,y)+
\sum_{m+b=1}^{\infty }\sum_{l=0}^m(
{\hat P}_{u^{b+m-l},w(u,\omega )^l}[a_{m-l+b,l}
(u,x(u,y),y)\circ $$ $$(I^{\otimes b}\otimes a^{\otimes (m-l)}
\otimes E^{\otimes l})])|_{u=t},$$ consequently,
$$\| X_1(t,y)-X_0(t,y) \| ^g \le $$ $$\sup_{m,l,b} (\| 
{\hat P}_{u^{b+m-l},w(u,\omega )^l}[a_{m-l+b,l}
(u,x(u,y),y)\circ (I^{\otimes b}\otimes a^{\otimes (m-l)}
\otimes E^{\otimes l})])|_{u=t} \| ^g .$$
Due to Condition $(ii)$ 
for each $\epsilon >0$ and $0<R_2<\infty $ there exists
$B_{\epsilon }\subset B_R$ such that 
$$K \sup_{m,l,b} (\| 
{\hat P}_{u^{b+m-l},w(u,\omega )^l}|_{B_{\epsilon }} 
[a_{m-l+b,l}(u,*,y)|_{(B_{\epsilon }\times B(L^q,0,R_2)\times V)}
\circ (I^{\otimes b}\otimes a^{\otimes (m-l)}
\otimes E^{\otimes l})]) \| ^g$$ 
$=:c<1.$
On the other hand, the partial difference quotient
has the continuous extension
$\bar \Phi ^1(X_{n+1}-X_n)(y;h;\zeta )$, that is expressible through
$\bar \Phi ^1$ of $a_{l,k}$, $a$ and $E$, and also through 
$a_{l,k}$, $a$ and $E$ themselves, where $y\in V$, $h\in Y$, 
$\zeta \in \bf K$ such that $y+\zeta h\in V$, since
analogous to $(X_{n+1}-X_n)$ estimates are true for 
$\bar \Phi ^1(X_{n+1}-X_n)$.
Therefore, there exists the unique solution on each
$B_{\epsilon }$ and it is of class $C^1$ by $y\in V$, 
since $\sup_{u,y}\max ( \| X_1(u,y)-X_0(u,y) \|_{L^q(\Omega ,H)}, 
\| \bar \Phi ^1(X_1(u,y)-X_0(u,y)) \|_{L^q(\Omega ,H)}<\infty $
and $\lim_{l\to \infty }c^lC=0$ for each $C>0$,
hence there exists $\lim_{n \to \infty }X_n(t,y)=X(t,y)=\xi (t,\omega ,y)
|_{B_{\epsilon }}$, where
$C:=M\sup_{u\in B_{\epsilon },y\in V}\max (\| X_1(u,y)-X_0(u,y) 
\| ^q_{L^q(\Omega ,H)}, 
\| \bar \Phi ^1 (X_1(u,y)-X_0(u,y)) \| ^q_{L^q(\Omega ,H)}
\le (c+1)K<\infty ,$ 
here $B_{\epsilon }$ is an arbitrary ball of radius 
$\epsilon $ in $B_R$, $t\in B_{\epsilon }$. 
Therefore, $\xi (t,\omega ,y)\in  C^1(V,L^q(\Omega ,{\sf F},\lambda ;
C^0(B_R,H))$.
\par From Proposition 3.11 \cite{lunast2} it follows, that
the multiplicative operator functional $T(t,v;\omega ;y)$
is of class $C^1$ by the parameter $y\in V$ such that
$\xi (t,\omega ,y)=T(t,v;\omega ;y)\xi (v,\omega ,y)$
for each $t$ and $v \in T$.
\par Due to the existence and uniqueness of the solution
$\xi (t,\omega ,y)$ for each $y\in V$, there exists
the operator $U(y_2,y;\xi (t,\omega ,y)):=\xi (t,\omega ,y_2)$,
that may be nonlinear by $\xi .$ The variation of the family
of solutions $\{ \xi (t,\omega ,y): y \} $ corresponds to the 
differential $D_y\xi (t,\omega ,y)$. Since $\xi (t,\omega ,y)$ 
is of class $C^1$ by $y$, then $U(y_2,y;\xi (t,\omega ,y)$
is of class $C^1$ by $y$ and $y_2$. 
The operator $U(y_2,y;*)$ has the inverse, since
$U(y,y_2;U(y_2,y;\xi (t,\omega ,y)))=\xi (t,\omega ,y)$
for each $y_2$ and $y\in V$, $t\in T$ and $\omega \in \Omega $.
Therefore, $U^{-1}(y_2,y;*)$ is also of class
$C^1$ by $y_2$ and $y$. In view of Conditions $(iii,iv)$
and $\xi (t,\omega ,y_2)-\xi (t,\omega ,y)\in X_{0,s}(H)$.
On the other hand, either $\mu _{S,\gamma ,q}$
or $\mu _S$ is quasi-invariant relative to shifts
$z\in X_{0,q}(C^0_0(T,H))$ and $S=S_1\otimes S_2$, consequently,
the transition measure $P_y$ is quasi-invariant relative to 
shifts $z\in X_{0,s}(H)$. In view of Conditions $(ii-iv)$ 
$\partial U(y_2,y;\eta )/\partial \eta -I\in L_r(H)$
for each $y_2$ and $y\in V$, where  
$\eta \in \{ \xi (t,\omega ,y): y \} $,
either $r=q$ or $r=0$ respectively.
Since $\mu _S(C^0_0(T,H))=1$, then $P_y(H)=1$, hence
$U(y_2,y;*)$ is defined $P_y$-almost everywhere
on $H$ for each $y_2$ and $y\in V$.
Therefore, there exists $n$ such that for each
$j>n$ the mappings $V(j;x):=x+P_j(U^{-1}(x)-x)$
and $U(j;x):=x+P_j(U(x)-x)$ are invertible and
$\lim_j|det U'x(j;x)|=|det U'_x(x)|$ and
$\lim_j|det V'_x(j;x)|=1/|det U'_x(x)|$, where
$U(x):=U(y_2,y;x)$, $y_2$ and $y\in V$. 
\par In view of Theorem 3.28 \cite{lu6} for each
$y_2$ and $y\in V$ the transition measures $P_{y_2}$ 
and $P_y$ are equivalent.
\par {\bf 2.14. Theorem.} {\it Let Conditions $2.13.(i-iv)$ be satisfied
and let $\phi $ be a $C^1$-diffeomorphism of a subset $V$ clopen in $\bf K$
onto the unit ball $B({\bf K},0,1)$. 
Then
\par $(1)$ the transition measure $P_y$ corresponding to 
$\mu _{S,\gamma ,q}$ is pseudo-differentiable by the parameter $y=\phi (z)$
of order $b\in \bf C$ for each $Re (b)\ge 0$, where $z\in V$;
\par  $(2)$ $P_y$ corresponding to $\mu _S$ with $h_k$ such that
$\sum_k\delta _k<\infty $, where 
$\delta _k:=\sup_{x\in B({\bf K},0,1)} |PD_c(b,h_k(x))|$,
is pseudo-differentiable by the parameter $y=\phi (z)$ of order $b$
for each $b\ge 0$, moreover,
$P_y$ is pseudo-differentiable for each $b\in C$, when
each $f_k$ is locally constant, that is, $h_k=0$ for each $k\in \bf N$.}
\par {\bf Proof.} Up to a constant multiplier
the operator $PD_c(b,h(x))$ of \S 2.7 
coincides with the pseudo-differential operator
$D^b(h(x)Ch_{B({\bf K},0,1)}(x))$ from \S 9 \cite{vla3}, 
where $Ch_A$ is the characteristic 
function of the subset $A$ in $\bf K$.
If $\psi \in L^2({\bf K},w,{\bf C})$ and $b>0$, then due to the 
Cauchy-Schwarz-Bunyakovskii inequality there exists
$\int_{{\bf K}\setminus B({\bf K},x,1)}[\psi (x)-\psi (y)]
|x-y|^{-1-b}w(dy)$, where $w$ is the Haar nonnegative measure 
on $\bf K$.
Then $F[D^b(h(x))]=|x|^bF[h(x)],$
where $F(h)(x):=\int_{\bf K}h(y) \chi _1((x,y))w(dx)$
is the Fourier transform  \cite{vla3,roo}
(see also \S 3.6 \cite{lunast1}). 
In view of Theorem 7.4 \cite{vla3} the Fourier transform
$f\mapsto F[f]$ is the bijective continuous isomorphism
of $L^2({\bf K},w,{\bf C})$ onto itself such that
$f(x)=\lim_{r\to \infty }\int_{B({\bf K},0,r)}F[f](y)
\chi _1(-(y,x))w(dy)$ and $(f,g)=(F[f],F[g])$
for each $f, g \in L^2({\bf K},w,{\bf C}).$
If $F[\psi ](x)=C exp (-\beta |x|^q)\chi _{\gamma }(x)$,
then there exists $D^b\psi (x)$ for each $b\ge 0.$
In accordance with Example $4.3.9$ $\int_{\bf K}\chi _{\gamma }(x)w(dx)=0$ 
for each $\gamma \ne 0$.
In view of Example $4.3.10$ \cite{vla3} 
$\int_{\bf Q_p}|x|^{nq}\chi _1(yx)w(dx)=[1-p^{nq}][1-p^{-n(q+1)}]^{-1}
|y|^{-n(q+1)}$ for each $q\in \bf C$ with $Re (q)>0$
and $n\in \{ 1,2,3,... \} $. 
\par If $f$ is a locally constant function as in \S 2.13,
then $PD_c(b,f)$ exists for each $b\in \bf C$.
On the other hand, $PD_c(b,f+h_k)=PD_c(b,f)+PD_c(b,h_k)$.
\par Let $g$ be a continuously differentiable
function $g: {\bf R}\to \bf R$ such that
$\| g \| _{C^1({\bf R},{\bf R})}:=\sup_x |g(x)| +\sup_x|g'(x)| <\infty $,
that is $g\in C^1_b({\bf R},{\bf R}).$
If for $f: {\bf K}\to \bf R$ and $x\in \bf K$
there exists $[f(x)-f(y)]|x-y|^{-1-b}
\in L^1({\bf K},w,{\bf C})$ as the function by $y\in \bf K$,
then $\int_{\bf K}[g\circ f(x)-g\circ f(y)] |x-y|^{-1-b}w(dy)=
\int_{S(f,x)}[g\circ f(x)-g\circ f(y)] [f(x)-f(y)]^{-1}
[f(x)-f(y)]|x-y|^{-1-b}w(dy)$, where $S(f,x):=\{ y: y \in {\bf K}, 
f(x)\ne f(y) \} $, consequently, there exists
$PD(b,g\circ f)(x)$. 
\par If instead of $g$ there exists $h\in C^1({\bf K},{\bf K})$
such that $\| h \| _{C^1({\bf K},{\bf K})}:=$ \\
$\max (\sup_x|h(x)|,
\sup_{x,y}| \bar \Phi ^1h(x;1;y)|)<\infty $, that is
$h\in C^1_b({\bf K},{\bf K})$, then
$\int_{\bf K}[f\circ h(x)-f\circ h(y)]|x-y|^{-1-b}w(dy)
=\int_{S(h,x)}[f\circ h(x)-f\circ h(y)]|h(x)-h(y)|^{-1-b}
|h(x)-h(y)|^{1+b} |x-y|^{-1-b}w(dy)$ exists,
hence there exists $PD(b,f\circ h)(x)$.
Analogous two statements are true for the operator $PD_c$ instead of $PD$.
\par In view of Equation 9.(1.5) \cite{vla3} $D^{\alpha }D^{\beta }\psi
=D^{\beta }D^{\alpha }\psi =D^{\alpha +\beta }\psi $
for each $\alpha \ne -1$, $\beta \ne -1$ and $\alpha +\beta \ne -1$
for each $\psi \in \sf D'$ such that there exists
$D^{\alpha }\psi $, $D^{\beta }\psi $ and $D^{\alpha +\beta }\psi $,
where $\sf D'$ is the topologically dual space
to the space $\sf D$ of locally constant functions
$\phi : {\bf K}\to \bf R$. On the other hand, 
$\sf D$ is dense in $\sf D'$
in the weak topology (see \S 6 \cite{vla3}).
Evidently, $L^2\cap \sf D$ is dense in $L^2({\bf K},w,{\bf R})$ also.
The characteristic functional of the Gaussian measure
belongs to $\sf D'$ and is locally constant 
on ${\bf K}\setminus \{ 0 \} $. Due to \S \S 7.2 and 7.3
the Fourier transform is the linear topological isomorphism
of $\sf D$ on $\sf D$ and of $\sf D'$ on $\sf D'$.
Then $\mu ^g_{S,\gamma ,q}(dx)/w(dx)\in 
L^1({\bf K},w,{\bf R})\cap \sf D'$ for each $g\in C^0_0(T,H)^*$.
\par In view of Theorem 4.3 \cite{lu6} and 
using the Kakutani theorem as in \S 2.13  
we get the statements of this theorem, since 
the quasi-invariance factor $P_y(dx)/P_u(dx)$ is 
pseudo-differentiable as the function by $y$ of order $b$ 
for each fixed $u\in B({\bf K},0,1)$.
\par {\bf 2.15. Theorem.} {\it Let $G$ be either a loop group
or a diffeomorphism group defined as in \S \S 2.9 and 2.11,
then there exists a stochastic process $\xi (t,\omega )$
on $G$ which induces a quasi-invariant transition measure $P$ 
on $G$ relative to $G'$ and $P$ is
pseudo-differentiable of order $b$ for each $b\in \bf C$ such that
$Re (b)\ge 0$ relative to $G'$, where
a dense subgroup $G'$ is given in \S \S 2.10 and 2.11.}
\par {\bf Proof.} These topological groups also have structures 
of $C^{\infty }$-manifolds, which are infinite-dimensional over the local 
field $\bf K$, but they do not satisfy the 
Campbell-Hausdorff formula in any open local subgroup
\cite{lutmf99,lubp2}.
Their manifold structures and actions of $G'$ on $G$ will be sufficient
for the construction of desired measures. 
These separable Polish groups have embeddings as clopen subsets into
the corresponding tangent Banach spaces $Y'$ and $Y$ in accordance with
\cite{luum985} and \S \S 2.1-2.11, where $Y'$ is the dense subspace
of $Y$. As usually $TG=\bigcup_{x\in G}T_xG$ and $T_xG=(x,Y)$.
\par Let $G$ be a complete separable relative 
to its metric $\rho $ $C^{\infty }$-manifold
on a Banach space $Y$ over $\bf K$ such that it 
has an embedding into $Y$ as the clopen subset.
Let $\tau _{G}: TG\to G$ be a tangent bundle on $G$. 
It is trivial, since $TG=G\times Y$ for the considered here case.
Let $\theta : Z_{G}\to G$ be a trivial bundle on $G$
with the fibre $Z$ such that $Z_{G}=Z\times G$, then $L_{1,2}(
\theta , \tau _{G})$ be an operator bundle with a fibre
$L_{1,2}(Z,Y)$ (see \S 2.13 \cite{lunast1}).
Let $\Pi :=
\tau _{G}\oplus L_{1,2}(\theta ,\tau _{G})$ be a Whitney sum
of bundles $\tau $ and $L_{1,2}(\theta ,\tau _{G})$. 
\par Since $G$ is clopen in $Y$, the valuation group of $\bf K$ 
is discrete in $(0,\infty )$, then it has a clopen
disjoint covering by balls $B(Y,x_j,r_j)$. That is, the atlas $At(G)$
of $G$ has a refinement $At'(G)$ being a disjoint atlas.
\par On $Y$ consider the measure $\mu _{S,\gamma ,q}$ or $\mu _S$ 
as in \S 2.12. Then in view of Theorems 4.3 \cite{lunast1}
and 2.2 \cite{lunast2} there exists the stochastic process
$w(t,\omega )$ corresponding to $\mu _{S,\gamma ,q}$ 
or $\mu _S$ (see also Definitions 4.2 \cite{lunast1} 
and 3.2 \cite{lunast2}).
Suppose that $f$ and $h_k$ for each $k\in \bf N$
defining the measure $\mu _S$ satisfy the Conditions
of \S 2.12 and of Theorem 2.14.
\par Now let $G$ be a loop or a diffeomorphism group of 
the corresponding manifolds 
over the field $\bf K$. 
Consider for $G$ a field $\sf U$ with a principal part
$(a_{\eta },E_{\eta })$, where $a_{\eta }\in T_{\eta }G$ and
$E_{\eta } \in L_{1,2}(H,T_{\eta }G)$ and $ker (E_{\eta })= \{ 0 \} $, 
$\theta : H_G\to G$ is a trivial bundle 
with a Banach fiber $H$ and $H_G:=G\times H$,
$L_{1,2}(\theta ,\tau _{\eta })$
is an operator bundle with a fibre $L_{1,2}(H,T_{\eta }G)$
such  that $(a_{\eta },E_{\eta })$ satisfies
Conditions of Theorem 3.4 \cite{lunast2}.
For Equation 3.5.(i) \cite{lunast2} we take additionally
$(a_{l,k})_{\eta }$ for each $l, k$ 
satisfying conditions of Theorem 3.5 \cite{lunast2}.
To satisfy conditions of quasi-invariance and
pseudo-differentiability of transition measures theorems
we choose $a_{\eta }$, $E_{\eta }$ and $(a_{k,l})_{\eta }$
of class $C^1$ and satisfying Conditions 2.13.(iii,iv) 
by  $y:=\eta \in G'=:V$ for each $k, l$.
\par We can take initially $\mu _{I,s,q}$ or $\mu _I$
a cylindrical measure on a Banach space $X'$ such that
$T_{\eta }G'\subset X'\subset T_{\eta }G$. 
If $A_{\eta }$ is the $L_q$-operator or the $L_1$-operator
with $ker (A_{\eta })= \{ 0 \} $, then $A_{\eta }$ gives 
the $\sigma $-additive measure $\mu _{A_{\eta },A_{\eta }^*z,q}$
or $\mu _{A_{\eta }}$ in the completion $X'_{1,\eta }$ of $X'$
with respect to the norm $\| x \| _1:= \| A_{\eta }x \| $
(see also \S 2.12). 
\par There exists the 
solution $\xi (t,\omega ,\eta )=\xi _{\eta }(t,\omega )$ of 
stochastic antiderivational Equation 3.4.(i) or 3.5.(i) \cite{lunast2}.
When the embedding $\theta $ of $T_{\eta }G'$ into $T_{\eta }G$
is $\theta =\theta _1\theta _2$ with $\theta _1$ and $\theta _2$ 
of class $L_q$ for $\mu _{S,\gamma ,q}$ or of class $L_1$
for $\mu _S$, then there exists $A_{\eta }$ 
such that $\mu _{A_{\eta },A_{\eta }z,q}$ or $\mu _{A_{\eta }}$
is the quasi-invariant and pseudo-differentiable of order $b$
measure on $T_{\eta }G$
relative to shifts on vectors from $T_{\eta }G'$
(see Theorems 2.13 and 2.14).
Henceforth we impose such demand on $A_{\eta }$ for each $\eta \in G'$.
\par Consider left shifts $L_h: G\to G$
such that $L_h\eta :=h\circ \eta $. 
Let us take $a_e\in T_eG'$, $A_e\in L_{1,q}(T_eG',T_eG)$
or $A_e\in L_{1,1}(T_eG',T_eG)$ respectively, 
$(a_{k,l})_{\eta }\in L_{k+l}((T_eG)^{\otimes (k+l)};
T_eG)$ for each $k$ and each $l$,
where $H$, $TeG'$ and $T_eG$ in their own norm uniformities are 
isomorphic with $c_0(\omega _0,{\bf K})$.
Then we put $a_x=(DL_x)a_e$ and $A_x=(DL_x)\circ A_e$
for each $x\in G$, hence $a_x\in T_eG$ and 
$A_x\in L_{1,s}(H_x,(DL_x)T_eG)$, where $(DL_x)T_eG=T_xG$
and $T_eG'\subset T_eG$, $H_x:=(DL_x)T_eG'$, $s=q$ or $s=1$.
Operators $L_h$ are (strongly) $C^{\infty }$-differentiable
diffeomorphisms of $G$ such that $D_hL_h: T_{\eta }G\to 
T_{h\eta }G$ is correctly defined, since $D_hL_h=h_*$
is the differential of $h$. In view of the choice of $G'$ in $G$ each 
partial difference quotient ${\bar \Phi }^nL_h(X_1,...,X_n;
\zeta _1,...,\zeta _n)$ is of class $C^0$ and $D^nL_h$ is of class 
$L_{n+1,s}({TG'}^{\otimes n}\times G',TG)$ 
for each vector fields $X_1,...,X_n$ on $G'$, $\zeta _1,..,
\zeta _n\in \bf K$ with $\zeta _jp_2(X_j)+h\in G'$
and  $h\in G'$, since for each $0\le l\in \bf Z$ the embedding of $T^lG'$ 
into $T^lG$ is the product of two operators of the $L_q$-class
or the embedding is of the $L_1$-class, 
where $T^0G:=G$, $X=(x,X_x)\in T_xG$, $x\in G'$, $X_x\in Y'$,
$p_1(X)=x$, $p_2(X)=X_x$. Take a dense subgroup $G'$
from \S 2.10 or \S 2.11 correspondingly 
and consider left shifts $L_h$ for $h\in G'$.
\par The considered here groups $G$ are separable, 
hence the minimal $\sigma $-algebra generated by cylindrical 
subalgebras $f^{-1}({\sf B}_n)$, n=1,2,..., coincides with 
the $\sigma $-algebra $\sf B$ of Borel subsets of $G$, where 
$f: G\to \bf K^n$ are continuous functions, ${\sf B}_n$ 
is the Borel $\sigma $-algebra  of $\bf K^n$. Moreover, $G$ 
is the topological Radon space (see Theorem I.1.2 and 
Proposition I.1.7 \cite{dalf}).
Let 
\par $P(t_0,\psi ,t,W):=P( \{ \omega : \xi (t_0,\omega )=\psi ,
\xi (t,\omega )\in W \} )$ \\
be the transition probability of 
the stochastic process $\xi $ for $t\in T$, which is defined on a 
$\sigma $-algebra $\sf B$ of Borel subsets in $G$, $W\in \sf B$,
since each measure $\mu _{A_{\eta },A_{\eta }^*z,q}$
is defined on the $\sigma $-algebra
of Borel subsets of $T_{\eta }G$ (see above). 
On the other hand, $T(t,\tau ,\omega )gx=gT(t,\tau ,\omega )x$
is the stochastic evolution family of operators for each
$\tau \ne t\in T$.
There exists the transition measure $P(t_0,\psi , t,W)$ such that 
it is a $\sigma $-additive
quasi-invariant pseudo-differentiable of order $b$ relative to 
the action of $G'$ by the left shifts $L_h$ on $\mu $
measure on $G$, for example, $t_0=0$ and $\psi =e$
with the fixed $t_0\in T$ (see Definitions 2.6 and 2.7). 
\par {\bf 2.16. Note.} In \S 2.15 $G'$ is on the Banach space $Y'$
and $G$ on the Banach space $Y$ over $\bf K$ such that
$G'$ and $G$ are complete relative to their uniformities
${\sf U}_{G'}$ and ${\sf U}_G$. There are inclusions
$TG'=G'\times Y'\subset G\times Y'\subset G\times Y=TG$.
The completion of $TG'$ relative to the uniformity
${\sf U}_G\times {\sf U}_{Y'}$ produces the uniform space $G\times Y'$.
Therefore, each ${\sf U}_G\times {\sf U}_{Y'}$-uniformly continuous 
vector field
$X=(x,X_x)$ on $G'$ has the unique extension on $G$ such that
$X_x\in Y'$ for each $x\in G$ (see \S 8.3 \cite{eng}), 
where ${\sf U}_G|_{G'}\subset {\sf U}_{G'}$.
Thus the ${\sf U}_G\times {\sf U}_{Y'}$-$C^{\infty }$-vector field
$X$ on $G'$ has the ${\sf U}_G\times {\sf U}_{Y'}$-$C^{\infty }$-extension 
on $G$ and it provides the $1$-parameter group 
$\rho : {\bf K}\times G\to G$
of $C^{\infty }$-diffeomorphisms of $G$ generated by a 
${\sf U}_G\times {\sf U}_{Y'}$-$C^{\infty }$-vector field 
$X_{\rho }$ on $G'$ \cite{lud,lutmf99}, that is,
$(\partial \rho (v,x)/\partial v)|_{v=0}=X_{\rho }(x)$ 
for each $x\in G$, where $v\in \bf K$, $X_{\rho }(x)\in G\times Y'$. 
In view of \S 2.15 the transition measure $P$ is quasi-invariant and 
pseudo-differentiable of order $b$ relative to the $1$-parameter 
group $\rho $. 
\par This approach is also applicable to the case of two Polish manifolds
$G'$ and $G$ of class $C^{\infty }$ on $Y'$ and $Y$ over $\bf K$.
The quasi-invariance and pseudo-differentiability
of the measure $P$ on $G$ relative to the $1$-parameter group
$\rho $ (by the definition) means such properties of $P$ relative to the
${\sf U}_G\times {\sf U}_{Y'}$-$C^{\infty }$-vector field 
$X$ on $G'$.
\par Evidently, considering different $(a,E)$ and $\{ a_{k,l}: k,l \}$
we see that there exist
${\sf c}=card ({\bf R})$ nonequivalent stochastic (in particular,
Wiener) processes on $G$ and $\sf c$ orthogonal quasi-invariant
pseudo-differentiable of order $b\in \bf C$ with $Re (b)>0$
measures on $G$ relative to $G'$.
\par If $M$ is compact, then in the case of the diffeomorphism group
its dense subgroup $G'$ can be chosen such that
$G'\supset Diff(t',M)$ for
$dim_{\bf K}M=n\in {\bf N}$ and $t'=t+s$ for $0\le t\in \bf R$,
$s>nv$, $v=dim_{\bf Q_p}({\bf K})$.
Analogously can be considered the manifold
$M\subset B({\bf K^n},0,r)$ and 
the group $G:=Diff(an_r,M)$ of analytic diffeomorphisms $f: M\to M$
having analytic extensions on $B({\bf K},0,r)$ with the corresponding
norm topology, where $r>0$ and $r<\infty $.
Then there exists the stochastic process  $\xi $ on $T_eG$ such that
it generates the transition measure $P$ on $T_eG$, its restriction
on the clopen subset $G$ embedded into $T_eG$ produces the quasi-invariant 
and pseudo-differentiable of each order $b\in \bf C$ with $Re (b)\ge 0$
measure $P|_G$ relative to the dense subgroup $G':=
Diff(an_R,M)$ for $R>r>0$, since the embedding 
$T_eG'$ into $T_eG$ is of class $L_1$ (see also \S 2.17).
\par {\bf 2.17.} {\bf Theorem.} {\it Let $G$ be a separable Banach-Lie
group over a local field $\bf K$. Then there exists a probability
quasi-invariant and pseudo-differentiable of each order $b\in \bf C$ with 
$Re (b)>0$ transition measure $P$ on $G$ relative to a dense subgroup
$G'$ such that $P$ is associated with a non-Archimedean stochastic process.}
\par {\bf Proof.} We consider two cases: $(I)$ $G$ satisfies locally the 
Campbell-Hausdorff formula; $(II)$ $G$ does not satisfy it in any 
neighbourhood of $e$ in $G$. The first case permits to describe $G'$ 
more concretely. There exists the embedding of $G$ into $T_eG$
as the clopen subset, since $G$ is the Polish group. 
The second case can be considered quite analogously to
\S 2.15, where the dense subgroup $G'$ can be characterized by the condition
that the embedding of $T_eG'$ into $T_eG$ is $\theta =\theta _1\theta _2$
with $\theta _1$ and $\theta _2$ of class $L_q$ or $\theta $ of class $L_1$, 
where $s=q$ or $s=1$ for stochastic processes associated 
either with $\mu _{S,\gamma ,q}$ or $\mu _S$ respectively.
\par It remains to consider the first case.
For $G$ there exists a Banach-Lie algebra $\sf g$
and the exponential mapping $exp: V\to U$, where $V$ is a neigbourhood
of $0$ in $\sf g$ and $U$ is a neighbourhood of $e$ in $G$ such that
$exp(V)=U$, where $exp(X+Y)=exp(X)exp(Y)$ for 
commuting elements $X$ and $Y$ of $\sf g$, that is, $[X,Y]=0$,
$exp(X)Yexp(-X)=exp(ad\mbox{ }X)Y$, $exp(\lambda X)=\sum_{j=0}^{\infty }
\lambda ^jX^j/j!$, $V=B({\sf g},0,r)$ is a ball of radius $0<r<\infty $
in $\sf g$, $\lambda \in \bf K$, $\lambda X\in V$, ${\sf g}=T_eG$.
The radii of convergence of the exponential and Hausdorff series 
corresponding to $log(exp(X).exp(Y))$ are
positive such that for each $0<R<p^{1/(1-p)}$
to a ball $B({\sf g},0,R)$
there corresponds a clopen subgroup $G_1$ supplied with the 
Hausdorff function (see \S II.6 and \S II.8  \cite{boua}).
Therefore, the exponential mapping $exp$ supplies $G$ with the 
structure of the analytic manifold over $\bf K$. By $At(G)=\{
(U_j,\phi _j): j \in {\bf N} \} $ is denoted the analytic atlas
of $\bf N$, that is $\phi _j: U_j\to V_j$ are diffeomorphisms
of $U_j$ onto $V_j$, where $U_j$ and $V_j$ are clopen in $G$
and in $\sf g$ respectively, connecting mappings 
$\phi _j\circ \phi _i^{-1}$ are analytic on $\phi _i(U_i\cap U_j)
\subset \sf g$. Therefore, the exponential mapping provides
$G$ with the covariant derivation $\nabla $ and a bilinear
tensor $\Gamma $ such that
$\nabla _XY=\mbox{ }_L\nabla _XY-\mbox{ }_LT(X,Y)/2$,
where the left-invariant derivation on $G$ is defined by
$\mbox{ }_L\nabla _X\tilde Y=0$ for an arbitrary left-invariant 
vector field $\tilde Y$ and all vector fields $X$ on $G$,
a vector field $\tilde Y$ is called left-invariant if $TL_g{\tilde Y}(h)=
{\tilde Y}(gh)$, $L_gh:=gh$ for each $g, h\in G$, $TL_g$ 
is the tangent mapping of $L_g$, $\nabla _XY_u=DY_u.X_u+
\Gamma _u(X_u,Y_u)$. For such $\nabla $
the torsion tensor is zero (see \S 1.7 \cite{kling},
\cite{freput} and \S 14.7 \cite{vla3}). It defines the rigid analytic 
geometry and the corresponding atlas on $G$.
Nevertheless $At(G)$ has the refinement $At'(G)$ 
such that charts of $At'(G)$ compose the disjoint covering of $G$.
\par Let $a_x$ be a an analytic vector field and $A_x$ be an anlytic
operator field on $G$ such that $A_x$ is 
an injective compact operator of class
$L_s$ for each $x\in G$, since $\sf g$
is of separable type over a spherically complete field
$\bf K$ and hence isomorphic with $c_0(\omega _0,{\bf K})$
(see Chapter 5 in \cite{roo}), where $s=q$ or $s=1$. 
Let $w_x(t,\omega )$ be a non-Archimedean stochastic 
(or, in particular, Wiener) process in $T_xG$ 
such that $a_xt+A_xw_x(t)\in T_xG$,
since the space $C^0_0$ is isomorphic with $c_0$.
For a ball $B_R:=B({\bf K},0,R)$ in $\bf K$ for $0<R<\infty $
let $B({\bf K},t_j,r)$ be a disjoint paving for sufficiently 
small $0<r<\infty $ for which $\xi ^q_x(t)=exp_{\xi ^q_{x,k}}
\{ a_{\xi ^q_{x,k}}(t-t_k)+A_{\xi ^q_{x,k}}[w_{\xi ^q_{x,k}}(t)
-w_{\xi ^q_{x,k}}(t_k) ] \} $ is defined, where $\xi ^q_{x,k}=
\xi ^q_x(t_k)$ for $k=0,1,...,n$, $\xi ^q_x(0)=x$, $q$ denotes the
partition of $B_R$ into $B({\bf K},t_j,r)$. Then there
exists the process $\xi =\lim_q\xi ^q(t)$ which is by our definition
a solution of the following stochastic equation:
$$(i)\mbox{ }d\xi (t,\omega )=exp_{\xi (t,\omega )} \{ a_{\xi (t,\omega )}dt+
A_{\xi (t,\omega )}dw(t,\omega ) \} $$
for $t\in B_R$. 
A function $f(t,x)$ such that 
$f(t,\xi ):=ln_{\xi (t,\omega )} \xi (t,\omega )$
satisfies the condition of Theorem 4.8 \cite{lunast1}
on the corresponding domain $W$, 
where $(t,x)\in W\subset T\times H$.
In view of Theorem 4.8 \cite{lunast1} after coordinate mapping
of a chart $(U,\phi )$ this equation takes the following form
on $\sf g$:
\par $$(ii)\quad \phi (\xi (t,\omega ))=\phi (\xi (t_0,\omega ))
+({\hat P}_ua^{\phi }_{\xi (u)})|_{u=t}
+({\hat P}_{w^{\phi }_{\xi (u)}(u,\omega )}E)|_{u=t},$$
$$-\sum_{m=2}^{\infty }(m!)^{-1}\sum_{l=0}^m{m\choose l}(
{\hat P}_{u^{m-l},w^{\phi }_{\xi (u)}
(u,\omega )^l}[(\partial ^{m-2}
\Gamma ^{\phi }_{\phi (\xi (u))}
/\partial x^{m-2})\circ (a^{\phi }_{\xi (u)}(u)^{\otimes (m-l)}
\otimes E^{\otimes l})])|_{u=t},$$
where $E=A^{\phi }_{\xi (u,\omega )}$,
$a^{\phi }=(\partial {\phi }_x/\partial x)a_x$, 
$A^{\phi }_x=(\partial {\phi }_x/\partial x)A_x
(\partial {\phi }^{-1}_x/\partial x)$, 
since $h^{\phi }=(\partial g^{\phi }/\partial x)f^{\phi }+\Gamma ^{\phi }_{
\phi (x)}(f^{\phi },g^{\phi })$ for $h=\nabla _fg$,
$f^{\phi }=(\partial {\phi }/\partial x)f$, $g^{\phi }=
(\partial {\phi }/\partial x)g$,
$h^{\phi }=(\partial {\phi }/\partial x)h$, $\Gamma ^{\phi }$ is 
a bilinear operator of Christoffel in $\sf g$,
which has the transformation property 
$D(\psi \circ \phi ^{-1}).\Gamma ^{\phi }_{\phi (x)}
=D^2(\psi \circ \phi ^{-1})+\Gamma ^{\psi }_{\psi (x)}\circ
(D(\psi \circ \phi ^{-1})\times D(\psi \circ \phi ^{-1}))$
such that $\nabla _XY_{\phi }=DY_{\phi }.X_{\phi }+
\Gamma ^{\phi }_{\phi (x)}(X_{\phi },Y_{\phi })$,
$\Gamma ^{\phi }_{\phi (x)}$ denotes $\Gamma $ for the chart
$(U,\phi )$, $\psi $ corresponds to another chart
$(V,\psi )$ such that $U\cap V\ne \emptyset ,$
$f,$ $g$ and $h$ are vector fields,
since $[\partial (\psi \circ \phi ^{-1})/\partial t]=0$, that is 
Corollary 4.7 \cite{lunast1} is applicable instead of Theorem 4.8 
\cite{lunast1} because $f$ corresponds to $(\psi \circ \phi ^{-1})$
(see \S 1.5 \cite{kling} and \cite{boum}). 
Since $a_x$ and $A_x$ are analytic, then 
$a$ and $E$ satisfy conditions of Theorem 3.4 \cite{lunast2}.
\par The function $\Gamma $ is analytic on the corresponding domain.
On the other hand, $\sf g$ is isomorphic with 
$c_0(\omega _0,{\bf K})$ as the Banach space.
If $Z$ is the center of $\sf g$, then $ad: {\sf g}/Z\to
gl(c_0(\omega _0,{\bf K}))$ is the injective representation,
where $gl(c_0)$ denotes the general linear algebra on $c_0$,
$ad(x)y:=[x,y]$ for each $x, y\in \sf g$.
Since $Z$ is commutative it aslo has injective representation
in $gl(c_0)$, consequently, $\sf g$ has embedding into
$gl(c_0(\omega _0,{\bf K})$, since $c_0\oplus c_0$ is isomorphic 
with $c_0$. Therefore, each $x\in \sf g$ can be written in the form
$x=\sum_{i,j}x^{i,j}X_{i,j}$, where $\{ X_{i,j}: i,j \in {\bf N} \} $ 
is the orthonormal basis of $\sf g$ as the Banach space, 
$x^{i,j}\in \bf K$, $\lim_{i+j\to \infty }x^{i,j}=0$, 
consequently, $\sf g$ has an embedding into 
$L_0(c_0(\omega _0,{\bf K})$.
Then $\Gamma $ can be written in local coordinates
$x_{s(i,j)}:=x^{i,j}$, where $s: {\bf N^2}\to \bf N$
is a bijection for which $\lim_{i+j\to \infty }s(i,j)=\infty $,
$X_{i,j}=:q_{s(i,j)}$,
since $(x^n)^{j,k}=\sum_{l_1,...l_{n-1}}x^{j,l_1}x^{l_1,l_2}...
x^{l_{n-1},k}X_{j,k}$, when $X_{j,k}=(\delta _{a,j}\delta _{b,k}: a, b
\in {\bf N} ) $, $\psi \circ \phi ^{-1}(x)=\sum_sq_sa^s_mx^m$
with $a^s_m\in \bf K$ and $\lim_{s+|m|+Ord(m)\to \infty }a^s_m=0$,
since $exp(x)$ has a radius of convergence $0<{\tilde r}=p^{-1}$
for $char({\bf K})=0$ (see Theorem 25.6 \cite{sch1}),
where $m=(m_1,...,m_k)$, $k=Ord(m)$, $0\le m_1\in \bf Z$,...,
$0\le m_{k-1}\in \bf Z$, $0<m_k\in \bf Z$, $0\le k\in \bf Z$.
Evidently, there exists $0<r<\infty $ such that the series
for $\psi \circ \phi ^{-1}$ converges in $B(c_0,0,r)$
for $V\cap U\ne \emptyset $. Hence each $a_{m-l,l}
:=(\partial ^{m-2} \Gamma ^{\phi }_{\phi (x)} /\partial ^{m-2}x)/m!$
for $m\ge 2$ and $a_{1,0}=a_{0,1}=(\partial \phi /\partial x)$ 
satisfies Condition $3.5.(ii)$ \cite{lunast2}. 
Due to Theorem $3.5$ \cite{lunast2} there exists the unique solution of
Equation $3.5.(i)$. 
Consider $G'$ corresponding to $\sf g'$ such that
the embedding $\theta $ of $\sf g'$ into $\sf g$ is of class $L_1$ for 
$\mu _S$ or $\theta =\theta _1\theta _2$, where $\theta _1$ and $\theta _2$
are of class $L_q$ for $\mu _{S,\gamma ,q}$. 
\par Let $T\in L_s({\sf g})$ or $T=T_1T_2$ with $T_1$ and $T_2
\in L_s({\sf g})$, where $s=1$ or $s=q$.  
Consider ${\sf h}_1:=T({\sf g})$, ${\sf h}_2:=sp_{\bf K}([{\sf h}_1,
{\sf g}]\cup {\sf h}_1)$ and by induction
${\sf h}_{n+1}:=sp_{\bf K}([{\sf h}_n,
{\sf g}]\cup {\sf h}_n)$, then ${\sf h}_{n+1}\supset
{\sf h}_n$ and ${\sf h}_n$ is the subalgebra in $\sf g$ for 
each $n\in \bf N$. In view of Proposition 2.12.2 the space 
$L_s({\sf g})$ is the ideal in $L({\sf g}).$
Therefore, ${\sf h}:=\bigcup_n{\sf h}_n$ is the ideal in $\sf g$
due to the anticommutativity and the Jacobi identity.
Since $\bf K$ is spherically complete there exists ${\sf h}_{n+1}\ominus 
{\sf h}_n=:{\sf t}_{n+1}$ for each $n\in \bf N$ and ${\sf t}_1:=
{\sf h}_1$ such that ${\sf t}_n$ is the $\bf K$-linear subspace of $\sf g$
(\cite{roo}). By ${\sf c}_0({\sf g}, \{ {\sf t}_n: n \} )=:\sf y$
is denoted the completion in $\sf g$ of vectors $z$ such that
$z=\sum_nz_n$ with $z_n\in {\sf t}_n$ for each $n$ and 
$\lim_{n\to \infty }z_n=0$. Evidently $\sf y$ is the proper ideal in 
$\sf g$ such that ${\sf h}\subset \sf y$, since
$\sf g$ is infinite dimensional over $\bf K$.
Then the embedding $\theta $ of $\sf y$ into $\sf g$ is either of class
$L_1$ or $\theta =\theta _1\theta _2$ such that
$\theta _1$ and $\theta _2$ belong to $L_q$. 
\par Due to this let $a$ and $A$ be such that
$[a_x,{\sf y}]\subset \sf y$ and $[A_x, ad ({\sf y})]
\subset ad ({\sf y})$ for each $x\in G$, where
$ad (x)g:=[x,g]$ for each $x, g \in \sf g$, that is, 
$ad (x) \in L({\sf g})$. If $g\in {\sf y}\cap V$, then 
$exp (ad (g)) - I$ is either of the class $L_1$ or
the product of two operators each of which is of the class $L_q$.
\par There exists a countable family $(g_j,W_j): j\in {\bf N})$ 
of elements $g_j\in G\setminus W$ for each $j>1$ and clopen subsets
$e\in W_j\subset W$ such that $g_1=e$, $W_1=W$ and $\{ g_jW_j: j \} $ is a locally finite covering
of $G$, since $G$ is separable and ultrametrizable (see \S 5.3 \cite{eng}).
If $P$ is a quasi-invariant and pseudo-differentiable of order $b$
measure on a clopen subgroup 
$W$ relative to a dense subgroup $W'$, then $P(S):=
(\sum_jP((g_j^{-1}S)\cap W_j)2^{-j}) (\sum_jP(W_j)2^{-j})^{-1}$
for each Borel subset $S$ in $G$ is 
quasi-invariant and pseudo-differentiable of order $b$
measure on $G$ relative to the dense subgroup
$G':=\bigcup_jg_j(W_j\cap W')$.
The group $G$ is totally disconnected and 
is left-invariantly ultrametrizable (see \S 8 and Theorem 5.5
\cite{hewr} and \S 6.2 \cite{eng}), consequently, in each neighbourhood
of $e$ there exists a clopen subgroup in $G$.
Then conditions of Theorems 2.13 and 2.14 are satisfied.
Therefore, analogously to \S 2.15
there are $S$, $\gamma $ and the stochastic process corresponding to 
$\mu _{S,\gamma ,q}$ or $\mu _S$ such that the transition measure
$P$ is quasi-invariant and pseudo-differentiable  relative to $G'$.
\par {\bf 2.18.} Theorem 2.15 gives the subgroup $G'$ concretely
for the given group $G$, but Theorem 2.17 describes concretely
$G'$ only for the case of $G$ satisfying the Campbell-Hausdorff formula.
For a Banach-Lie group not satisfying locally the Campbell-Hausdorff 
formula Theorem 2.17 gives only the existence of $G'$.
\par These transition measures $P=:\nu $ on $G$
induce strongly continuous unitary regular representations
of $G'$ given by the following formula:
$T_h^{\nu }f(g):=(\nu ^h(dg)/\nu (dg))^{1/2}f(h^{-1}g)$
for $f\in L^2(G,\nu ,{\bf C})=:H$, $T_h^{\nu }\in U(H)$,
$U(H)$ denotes the unitary group of the Hilbert space $H$. 
For the strong continuity of $T_h^{\nu }$
the continuity of the mapping $G'\ni h\mapsto \rho _{\nu }(h,g)\in 
L^1(G,\nu ,{\bf C})$ and that $\nu $ is the Borel measure
are sufficient, where $g\in G$, since $G$ is the Polish space and 
hence the Radon space (see Theorem I.1.2 \cite{dalf}). 
On the other hand, the continuity of $\rho _{\nu }(h,g)$ by $h$
from the Polish group $G'$ into $L^1(G,\nu ,{\bf C})$ follows from
$\rho _{\nu }(h,g)\in L^1(G,\nu ,{\bf C})$  for each 
$h\in G'$ and that $G'$ is the topological subgroup of $G$
(see \cite{fidal}).
\par Then analogously to \S 2.15 there can be constructed 
quasi-invariant and pseudo-differentiable measures on the manifold
$M$ relative to the action of the diffeomorphism group 
$G_M$ such that $G'\subset G_M$. Then Poisson measures on 
configuration spaces associated with either $G$ or $M$ can be
constructed  \cite{lupmum}.
There exists the stochastic process corresponding to
$\mu _{S,\gamma ,q}$ with the certain choice of $a$, $E$ and $a_{k,l}$
such that the regular representation is irreducible, for 
the stochastic process corresponding to $\mu _S$ it can be taken
the family of $\{ f_k: k \} $ and $a$, $E$ and $a_{k,l}$ 
such that the regular representation 
is irreducible. 
\par More generally it is possible to consider
instead of the group $G$ a Polish topological space
$X$ on which $G'$ acts jointly continuously: $\phi : (G'\times X)\ni
(h,x)\mapsto hx=:\phi (h,x)\in X$, $\phi (e,x)=x$
for each $x\in X$, $\phi (v,\phi (h,x))=\phi (vh,x)$
for each $v$ and $h\in G'$ and each $x\in X$.
If $\phi $ is the Borel function, then
it is jointly continuous \cite{fidal}.
From \S II.3.2 \cite{lubp2} (see aslo
\cite{lutmf99,lurim2,lupmum}) there is the following result.
\par {\bf Theorem.} {\it Let $X$ be a Polish topological space 
with a $\sigma $-additive $\sigma $-finite nonnegative nonzero 
ergodic Borel measure
$\nu $ quasi-invariant relative to a Polish topological group
$G'$ acting on $X$ by the Borel function $\phi $. If 
\par $(i)$ $sp_{\bf C} \{\psi |
\quad \psi (g):=(\nu ^h(dg)/\nu (dg))^{1/2}, h\in G' \} $
is dense in $H$ and
\par $(ii)$ for each $f_{1,j}$ and $f_{2,j}$ in $H$, $j=1,...,n,$ 
$n\in \bf N$ and each $\epsilon >0$ there exists $h\in G'$ such that
$|(T_hf_{1,j},f_{2,j})| \le \epsilon |(f_{1,j},f_{2,j})|$,
when $|(f_{1,j},f_{2,j})|>0$.
Then the regular representation $T: G'\to U(H)$ is irreducible.}
\par There can be used pseudo-differentiable measures of 
order $l$ either for each $l\in \bf N$ or $-l\in \bf N$,
that is used for the verification of Condition $(i)$.
Transition measures corresponding to stochastic processes
that are quasi-invariant and pseudo-differentiable of each order
$b\in \bf C$ with $Re (b)\le 0$ can be constructed analogously
starting with the corresponding measures $\mu _S$.
To satisfy the conditions of this theorem,
for example, in \S 2.15 it can be taken $a=0$, $E$ nondegenerate
independent from $t$ and each $a_{k,l}=0$ besides $a_{0,1}=1$; in \S 2.17
it can be taken $a=0$, $E$ nondegenerate independent from $t$
and $a_{k,l}$ defined by the exponential mapping for $G$.
\par In view of Proposition II.1 \cite{neeb} for the separable 
Hilbert space $H$ the unitary group endowed with the strong 
operator topology $U(H)_s$ is the Polish group.
Let $U(H)_n$ be the unitary group with the metric
induced by the operator norm. In view of the Pickrell's theorem 
(see \S II.2 \cite{neeb}): if $\pi : U(H)_n\to U(V)_s$
is a continuous representation of $U(H)_n$ on  
the separable Hilbert space $V$, then $\pi $ is also
continuous as a homomorphism from $U(H)_s$ into
$U(V)_s$. Therefore, if $T: G'\to U(H)_s$ is
a continous representation, then there are new representations
$\pi \circ T: G'\to U(V)_s$. On the other hand, the 
unitary representation theory of $U(H)_n$ is the same as that of
$U_{\infty }(H):=U(H)\cap (1+L_0(H))$, since the group $U_{\infty }(H)$
is dense in $U(H)_s$.
\par Two theorems about induced representations of the dense
subgroups $G'$ were 
proved in \cite{luglgdg}, which are also applicable to the considered 
here cases.
\section{Stochastic antiderivational equations and measures
on a loop monoid and a path space.}
\par {\bf 3.1. Theorem.} {\it On the monoid 
$G=\Omega _{\xi }(M,N)$ from \S 2.3.1 and each $b_0\in \bf C$
with $Re (b_0)\ge 0$ there exists a stochastic process $\eta (t,\omega )$
on $G$ such that the transition measure $P$
is quasi-invariant and pseudo-differentiable of each order 
$b\in \bf C$ with $Re (b)\ge Re (b_0)$ 
relative to the dense submonoid $G':=\Omega _{\xi }^{ \{ k \} }(M,N)$ 
from \S 2.8 (with $c>0$ and $c'>0$).}
\par {\bf Proof.} In view of Lemma I.2.17 \cite{lubp2}
it is sufficient to consider the case of $M$ with the 
finite atlas $At'(M)$. The rest of the proof is quite analogous 
to that of Theorem 2.15 using the definitions of the 
quasi-invariance and the pseudo-differentiability for 
semigroups from \S 2.7.
\par {\bf 3.2. Definition and Note.} 
In view of \S 2.5 each space 
$N^{\xi }$ has the additive group structure,
when $N=B(Y,0,R)$, $0<R\le \infty $.
\par Therefore, the factorization by the
equivalence relation $K_{\xi }\times id$ produce the monoid of paths
$C^{\theta }_0(\xi ,\bar M\to N)/(K_{\xi }\times id)=:
{\sf S}_{\xi }(M,N)$ in which compositions are defined not for all
elements, where $y_1idy_2$ if and only if $y_1=y_2\in N$. 
There exists a composition $f_1f_2=(g_1g_2,y)$ if and only if
$y_1=y_2=y$, where $f_i=(g_i,y_i)$, $g_i\in \Omega _{\xi }(M,N)$
and $y_i\in N^{\xi }$, $i\in \{ 1,2 \} $.
The latter
semigroup has elements $e_y$ such that $f=e_y\circ f=f\circ e_y$ 
for each $f$, when their composition is defined, where $y\in N^{\xi }$,
$f=(g,y)$, $g\in \Omega _{\xi }(M,N)$, $e_y=(e,y)$.
If $N^{\xi }$ is a monoid, then ${\sf S}_{\xi }(M,N)$
can be supplied with the structure of a direct product
of two monoids.
Therefore, $P_{\xi }(M,N):=L_{\xi }(M,N)\times N^{\xi }$ is called
the path group.
\par {\bf 3.3. Theorem.} {\it On the path group $G=P_{\xi }(M,N)$ from
\S 3.2, when $N=B(Y,0,R)$ and $N^{\xi }$ is supplied with the 
additive group structure, and each $b_0\in \bf C$ with $Re (b_0)\ge 0$
there exists a stochastic process $\eta (t,\omega )$
for which a transition measure $P$ is quasi-invariant
and pseudo-differentiable of each order 
$b\in \bf C$ with $Re (b)\ge Re (b_0)$ relative to a dense
subgroup $G'$.}
\par {\bf Proof.} Since $P_{\xi }(M,N)=
L_{\xi }(M,N)\times N^{\xi }$, it is sufficient to construct 
two stochastic processes on $L_{\xi }(M,N)$ and $N^{\xi }$
and to consider transition measures for them.
In view of Theorems 2.15 and 2.17 the desired processes 
and transition measures for them exist.
\par {\bf 3.4. Definition.} Let the topology of $\Omega _{\xi }(M,N)$
be defined relative to countable $At(M)$.
If $F$ is the free Abelian group corresponding to
$\Omega _{\xi }(M,N)$, then there exists a set $\bar W$
generated by formal finite linear combinations over $\bf Z$ of elements
from $C^0_0(\xi ,(M,s_0)\to (N,y_0))$ and a continuous extension
$\bar K_{\xi }$ of $K_{\xi }$ onto $W_{\xi }(M,N)$ 
and a subset $\bar B$ of $\bar W$
generated by elements $[f+g]-[f]-[g]$ such that
$W_{\xi }(M,N)/\bar K_{\xi }$ is isomorphic with $L_{\xi }(M,N)$, where
$$W_{\xi }(M,N):=\bar W/\bar B,$$ 
$f$ and $g\in C^0_0(\xi ,(M,s_0)\to (N,y_0))$, $[f]$ is an element in
$\bar W$ corresponding to $f$, $\bar W$ is in a topology inherited from
the space $C^0_0(\xi ,(M,s_0)\to (N,y_0))^{\bf Z}$ in the Tychonoff product 
topology. We call $W_{\xi }(M,N)$ an $O$-group. Clearly the composition
in $C^0_0(\xi ,(M,s_0)\to (N,y_0))$ induces the composition in
$W_{\xi }(M,N)$. Then $W_{\xi }(M,N)$ is not the algebraic group,
but associative compositions are defined for its elements due to 
the homomorphism $\chi ^*$ given by Formulas 2.3.2.(5,6),
hence $W_{\xi }(M,N)$ is the monoid without the unit element.
\par Let $\mu _h(A):=\mu (h \circ A)$ for each $A\in Bf(W_{\xi }(M,N))$
and $h\in W_{\xi }(M,N)$, then as in \S 2.7 we get the definition
of quasi-invariant and pseudo-differentiable measures.
\par Let now $G':=W _{\xi }^{\{ k\} }
(M,N)$ be generated by $C_{0,\{ k\} }^0(\xi ,(M,s_0)\to (N,0))$ 
as in \S 2.8, then it is the dense $O$-subgroup in $W_{\xi }(M,N)$,
where $c>0$ and $c'>0$.
\par {\bf 3.5. Theorem.} {\it Let $G:=W_{\xi }(M,N)$ be the $O$-group
as in \S 3.4, $At(M)$ be finite and $b_0\in \bf C$ with $Re (b_0)\ge 0$.
Then there exists a stochastic process $\eta (t,\omega )$
on $G$ for which the transition 
measure $P$ is quasi-invariant and pseudo-differentiable 
of each order $b\in \bf C$ with
$Re (b)\ge Re (b_0)$ on $G$ relative to a dense $O$-subgroup $G'$.}
\par {\bf Proof.} The uniform space
$C^0_0(\xi ,M\to N)$ has the embedding as the clopen subset
into $C^0_0(\xi ,M\to Y)$. 
Here we can take $a\in TG'$ and $A\in L_{1,s}(\theta ,\tau )$ 
without relations with $DL_h$, where $s=q$ or $s=1$ respectively.
Then repeating the major parts of the proof of \S 2.15
without $L_h$ and so more simply, 
but using actions of vectors fields of $TG'$ by $\rho _X$ on $G$
we get the statement of this theorem, 
since $(D_X\rho _X)Y$ and $[(\nabla _X)^n(D_X\rho _X)]Y$
are products of two operators of
of class $L_{n+2,q}((TG')^{n+2},TG)$ and also 
of class $L_{n+2,1}((TG')^{n+2},TG)$
for each $C^{\infty }$-vector fields 
$X$ and $Y$ on $G'$ and each $n\in \bf N$.
In view of \S 2.16 there exists a stochastic process
$\eta (t,\omega )$ for which
the transition measure $P$ is quasi-invariant and 
pseudo-differentiable relative to each $1$-parameter diffeomorphism group
of $G'$ associated with a 
${\sf U}_G\times {\sf U}_{Y'}$-$C^{\infty }$-vector field on $G'$. 

\end{document}